\documentclass{article}
\usepackage{geometry}
\usepackage{amsmath,amsthm,amssymb,amsfonts}
\usepackage{mathtools}
\usepackage{physics}
\usepackage{booktabs}
\usepackage{graphicx}
\usepackage{tikz}
\usetikzlibrary{arrows, positioning, shapes}
\usepackage{xcolor}
\usepackage{url}
\usepackage{hyperref}      
\usepackage{bm}
\usepackage{tablefootnote}
\usepackage{authblk}
\usepackage{enumitem}

\theoremstyle{definition}
\newtheorem{theorem}{Theorem}[section]

\newtheorem{proposition}[theorem]{Proposition}
\newtheorem{lemma}[theorem]{Lemma}
\newtheorem{definition}[theorem]{Definition}
\newtheorem{example}[theorem]{Example}
\newtheorem{remark}[theorem]{Remark}

\newcommand{\FF}{\mathcal{F}} 

\newcommand{\define}[1]{{\bf \boldmath{#1}}}
\newcommand{\dist}{\mathrm{d}} 
\newcommand{\im}{\mathrm{im}} 
\newcommand{\birth}{\mathrm{birth}}
\newcommand{\death}{\mathrm{death}}

\newcommand{\N}{\mathbb{N}}        
\newcommand{\R}{\mathbb{R}}        
\newcommand{\Sp}{\mathbb{S}}       
\newcommand{\Hp}{\mathbb{H}}       
\newcommand{\Gr}{\mathbb{G}}       
\newcommand{\AG}{\mathbb{AG}}      

\newcommand{\CF}{\mathrm{CF}}      
\newcommand{\D}{\mathcal{D}}       
\newcommand{\radon}{\mathcal{R}}   

\newcommand{\PH}{\mathrm{PH}}      
\newcommand{\PD}{\mathrm{PD}}      

\newcommand{\PHT}{\mathrm{PHT}}                
\newcommand{\GPHT}{\PHT_{\mathbb{P},f}}        
\newcommand{\CPHT}{\PHT_{\mathbb{S}^{n-1}, h}} 
\newcommand{\DPHT}{\PHT_{\AG(m,n),\dist}}      

\newcommand{\HPHT}{\PHT_{\AG(n-1,n),\dist}}    
\newcommand{\TPHT}{\PHT_{\AG(1,n),\dist}}      
\newcommand{\RPHT}{\PHT_{\AG(0,n),\dist}}      
\newcommand{\ECT}{\mathrm{ECT}}

\author[1]{Adam Onus}
\author[2,3]{Nina Otter}
\author[2,3]{Renata Turke\v{s}}

\affil[1]{Queen Mary University of London, United Kingdom}
\affil[2]{DataShape, Inria-Saclay,  France}
\affil[3]{Laboratoire de Mathématiques d’Orsay,  Universit\'e Paris-Saclay, France}

\title{Shoving tubes through shapes gives a sufficient and efficient shape statistic}

\date{}

\begin{document}
\maketitle

\begin{abstract}
The Persistent Homology Transform ($\PHT$) was introduced in the field of Topological Data Analysis about 10 years ago, and has since been proven to be a very powerful descriptor of Euclidean shapes. The PHT consists of scanning a shape from all possible directions $v\in S^{n-1}$ and then computing the persistent homology of sublevel set filtrations of the respective height functions $h_v$; this results in a sufficient and continuous descriptor of Euclidean shapes. We introduce a generalisation of the PHT in which we consider arbitrary parameter spaces and sublevel sets with respect to any function. In particular, we study transforms, defined on the Grassmannian $\mathbb{A}\mathbb{G}(m,n)$ of affine subspaces of $\mathbb{R}^n$, that allow to scan a shape by probing it with all possible affine $m$-dimensional subspaces $P\subset \mathbb{R}^n$, for fixed dimension $m$, and by then computing persistent homology of sublevel set filtrations of the function $\mathrm{dist}(\cdot, P)$ encoding the distance from the flat $P$. We call such transforms ``distance-from-flat'' PHTs.
We show that these transforms are injective and continuous and that they provide computational advantages over the classical PHT. In particular, we show that it is enough to compute homology only in degrees up to $m-1$ to obtain injectivity; for $m=1$ this provides a very powerful and computationally advantageous tool for examining shapes, which in a previous work by a subset of the authors has proven to significantly outperform state-of-the-art neural networks for shape classification tasks.
\end{abstract}

\section{Introduction}
\label{sec:introduction}

\begin{figure}[h]
\centering
\begin{tabular}{c|c|c|c}
\toprule
Shape $X$ & $h_v$ & Height filtrations & $\PHT(X)$ \\
\midrule
\includegraphics[height = 0.1\linewidth]{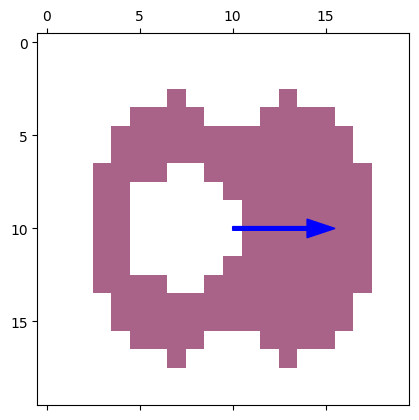} & 
\includegraphics[height = 0.1\linewidth]{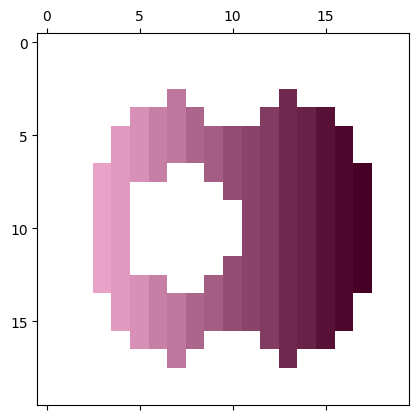} &
\includegraphics[height = 0.1\linewidth]{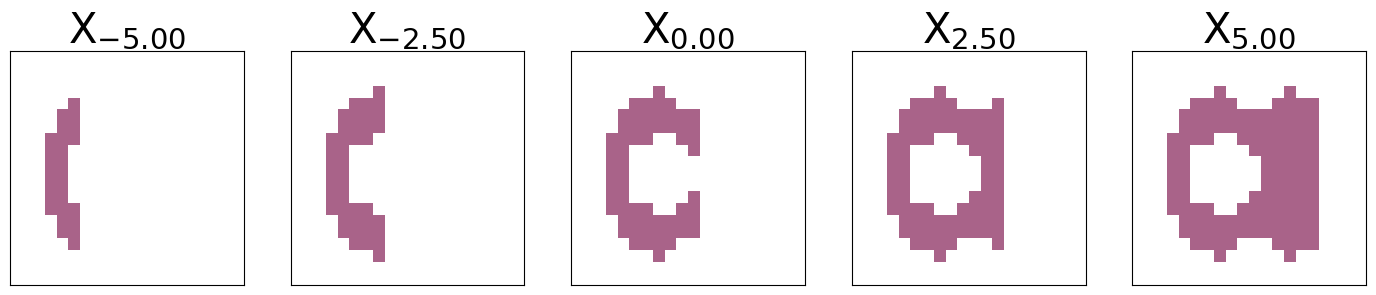} & 
\includegraphics[height = 0.1\linewidth]{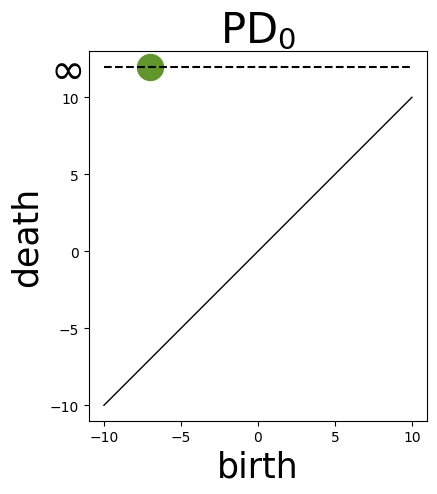} 
\includegraphics[height = 0.1\linewidth]{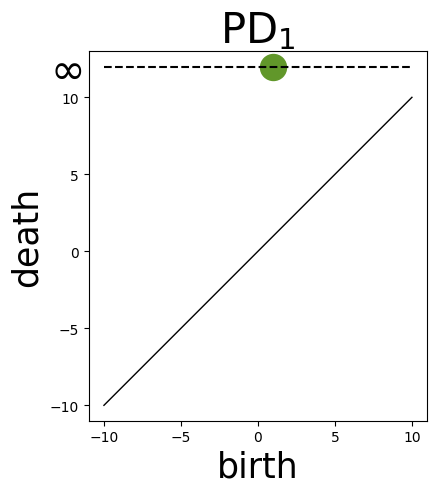} \\

\includegraphics[height = 0.1\linewidth]{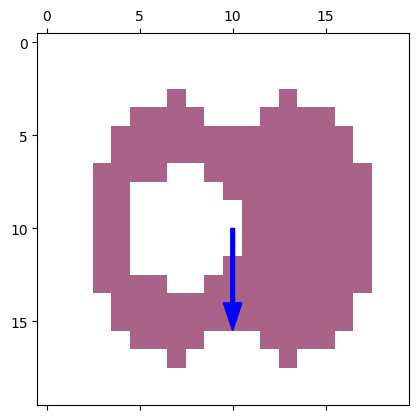} & 
\includegraphics[height = 0.1\linewidth]{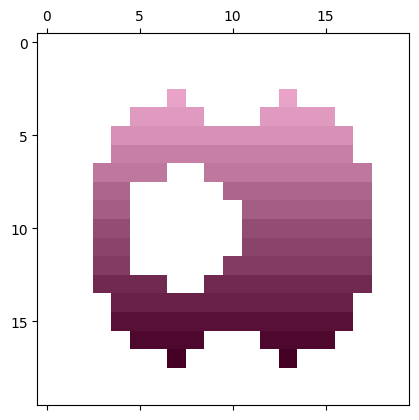} &
\includegraphics[height = 0.1\linewidth]{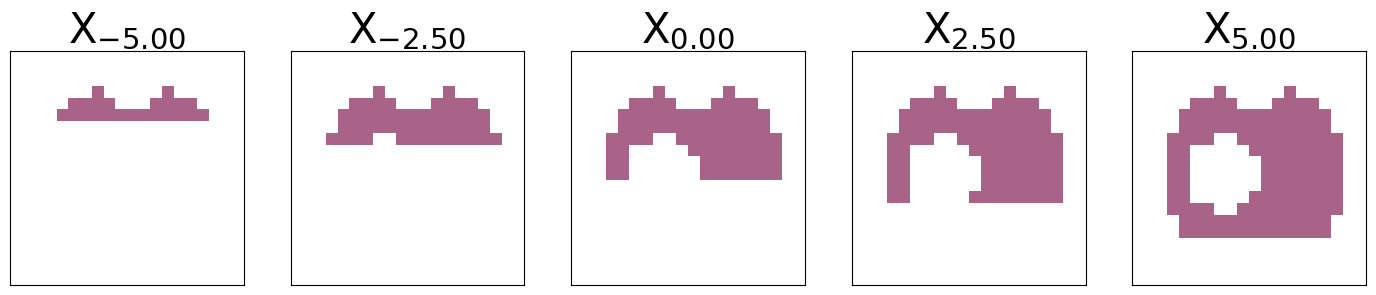} & 
\includegraphics[height = 0.1\linewidth]{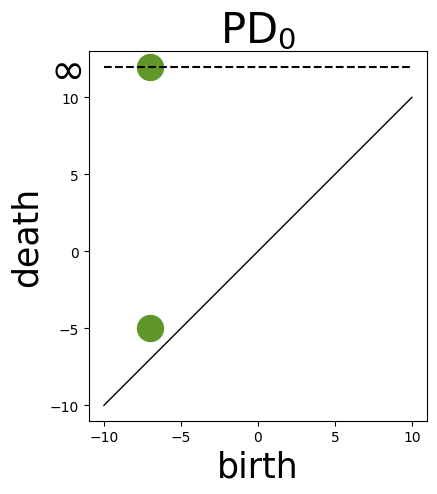} 
\includegraphics[height = 0.1\linewidth]{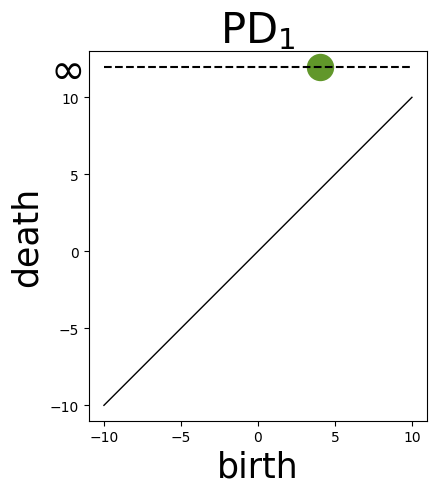} \\

... & ... & ... & ... \\
\bottomrule
\end{tabular}
\caption{For a shape $X$, the classical $\PHT(X)$ consists of $\PD_k(X, h_v)$ with respect to the height filtration $h_v$ for \emph{any} direction $v$, and for \emph{any} homological degree $k$. For the example image $X$, the figure illustrates $\PHT$ on the height filtration from two directions (in blue).}
\label{fig:pht}
\end{figure}

Shape classification is a difficult problem that plays a crucial role in understanding and recognising physical structures and objects, image processing, and computer vision \cite{lim2010shape}.  
Since topology is the branch of mathematics that studies shape, it has inspired the new field of study called Topological Data Analysis (TDA), which, as the name suggests,  aims to analyse data by studying its shape. 
The main tool in TDA, persistent homology ($\PH$), provides insights about the shape through information about $k$-dimensional cycles (connected components, loops or tunnels, voids and higher-dimensional voids). 
More precisely, given a nested family $\{\FF_r\}_{r \in \mathbb{R}}$ of subspaces that approximate the shape $X$ at any scale $r$ (a so-called ``filtration''), $\PH$ captures how the homological structures persist in the filtration as the parameter $r$ changes. The persistence of $k$-dimensional cycles is commonly summarised with a so-called ``persistence diagram'' $(\PD),$ a scatter plot of ``birth" and ``death" values of each cycle within the filtration. An example shape, filtration and resulting persistence diagrams are illustrated in any row of Figure~\ref{fig:pht}.
In this way, $\PH$ captures topological and geometric features about a shape \cite{tubular}, which are richer than the scalar-valued quantities that are often used in shape analysis and that struggle to describe the shape sufficiently well.

A remarkable result that highlights the effectiveness of $\PH$ in probing shapes was shown in \cite{TMB2014}. In that work, the persistent homology \emph{transform} $(\PHT)$, which we refer to as the \textbf{classical} $\PHT$, is introduced as a method that summarises persistent homology ($\PH$) information across  \emph{any} homological degree, with respect to the height filtration function, which scans $X\subset \mathbb{R}^{n}$ from \emph{any} given direction $v \in \mathbb{S}^{n-1}$ (Figure~\ref{fig:pht}). More precisely,
\begin{alignat*}{3}
\PHT(X) \colon & \mathbb{S}^{n-1} && \to \mathcal{D}^n \\ 
& v && \mapsto \bigg(\PD_0(X), \PD_1(X), \dots, \PD_{n-1}(X)\bigg) \, ,
\end{alignat*}
where $\PD_k(X)$ is the persistence diagram of $X$ in homological degree $k$, with respect to the height filtration function in direction $v$, i.e.,  $h_v(x)=x\cdot v$. 
The authors of \cite{TMB2014} consider finitely triangulisable subspaces of $\mathbb{R}^n$ and show that $\PHT(X)$ is continuous with respect to the Bottleneck and  Wasserstein distances on the space of persistence diagrams \cite[Lemma 2.1]{TMB2014}. Further, they  provide a constructive algorithmic proof of injectivity in dimension $n=2$ and $n=3$ \cite[Theorem 3.1]{TMB2014}. The injectivity result implies that a shape can be completely recovered from its $\PHT$.
Later,  an upper bound for the number of directions needed to determine a shape was given \cite[Theorem 7.14]{CST22}, as well as a  generalisation of the continuity and injectivity results to constructible sets, any Wasserstein distance $W_p$ and any dimension $n$ \cite[Lemma 4.8 and Theorem 4.16, respectively]{CST22}. A similar proof of injectivity was given at the same time in \cite{GLH2019}. Both injectivity proofs from \cite{CST22} and \cite{GLH2019} rely  on Schapira's inversion formula for integral transforms \cite{schapira1995}.

In this work, we extend beyond the classical $\PHT$ by allowing the following two crucial generalisations:
\begin{itemize}
\item arbitrary topological space $\mathbb{P}$ as the domain of the transform, and
\item arbitrary filtration functions $\{f_P\}_{P \in \mathbb{P}}$.
\end{itemize}
In particular, we focus on $\mathbb{P} = \AG(m, n)$, the  Grassmannian  of $m$-dimensional affine subspaces of $\mathbb{R}^n$ (which we call ``affine Grassmannian''), and $f_P = \dist_P(\cdot, P)$ the distance from a flat $P \in \AG(m, n),$ that we refer to as \textbf{distance-from-flat persistent homology transform}, and denote by $\DPHT$. For instance, if $m = 0$, $m=1$ or $m = n-1$, this corresponds to $\PH$ on respectively the distance from any point, line or hyperplane (which we refer to as radial-, tubular- and height-$\PHT$, see Table~\ref{tab:gpht}).

\begin{table}[h!]
\centering
\small
\begin{tabular}{lllll}
\toprule
Name & Notation & $\mathbb{P}$ & $\dim(\mathbb{P})$ & $k$ \\
\midrule
height $\PHT$ & $\HPHT$ & $\AG(n-1, n)$ = hyperplanes in $\R^n$ & n & 0, 1, \dots, n-2 \\
\dots & \dots & \dots & \dots & \dots \\
tubular $\PHT$ & $\TPHT$ & $\AG(1, n)$ = lines in $\R^n$ & 2(n-1) & 0 \\
radial $\PHT$ & $\RPHT$ & $\AG(0, n)$ = points in $\R^n$ & n & ``-1''\tablefootnote{Here we have that the Radon transform taking $X$ to $\chi(X\cap P)$ for all $0$-flats $P$ is its indicator function (see Lemma \ref{lemma:injectivity}).  Thus, we argue that one might intuitively think about this case as being homology in degree ``$-1$'': in the realm of $n$-categories, one would interpret a ($-1$)-category as truth values, see the periodic table of $n$-categories \cite[Section 2.5.2]{BS10}.} \\
\bottomrule
\end{tabular}
\caption{Relationship between distance-from-flat $\DPHT$ across different $m$ of the underlying affine Grassmannian $\AG(m, n)$. For low $m$, the dimension $\dim(\mathbb{P}) = \dim(\AG(m, n)) = (m+1)(n-m)$ \cite[Proposition 4.1]{lim2021grassmannian} of the parameter space remains low, in addition to less homological degrees $k$ that are sufficient to completely capture the shape.}
\label{tab:gpht}
\end{table}

\begin{figure}
\centering

\begin{tabular}{c|c|c|c}
\toprule
Shape $X$ & $d_P$ & Distance-from-flat filtrations & $\PHT(X)$ \\
\midrule

\includegraphics[height = 0.09\linewidth]{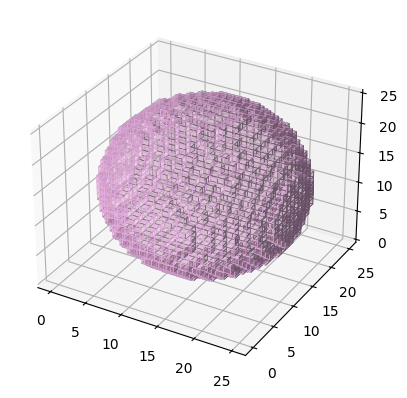} &
\includegraphics[height = 0.09\linewidth]{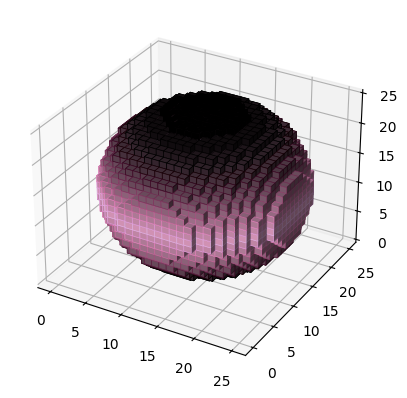} &
\includegraphics[height = 0.09\linewidth]{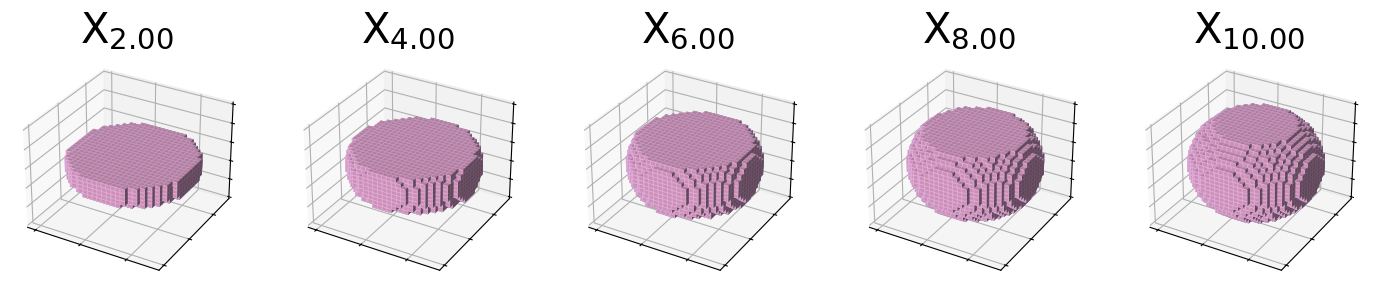} &
\includegraphics[height = 0.09\linewidth]{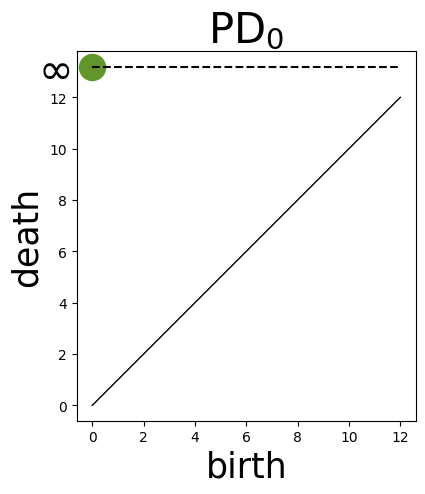} 
\includegraphics[height = 0.09\linewidth]{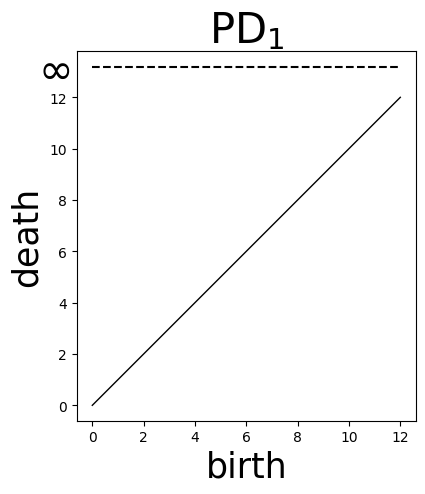} 
\includegraphics[height = 0.09\linewidth]{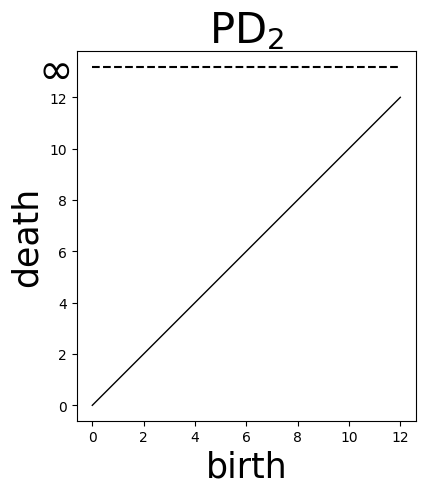} \\

\includegraphics[height = 0.09\linewidth]{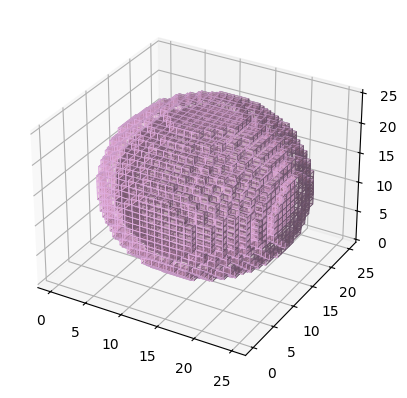} &
\includegraphics[height = 0.09\linewidth]{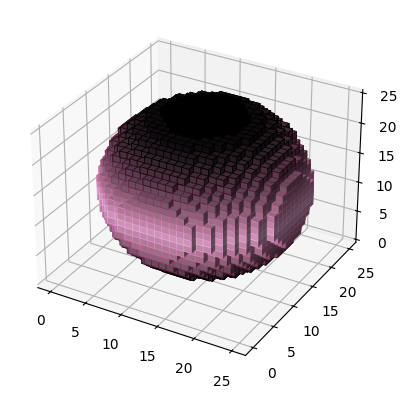} &
\includegraphics[height = 0.09\linewidth]{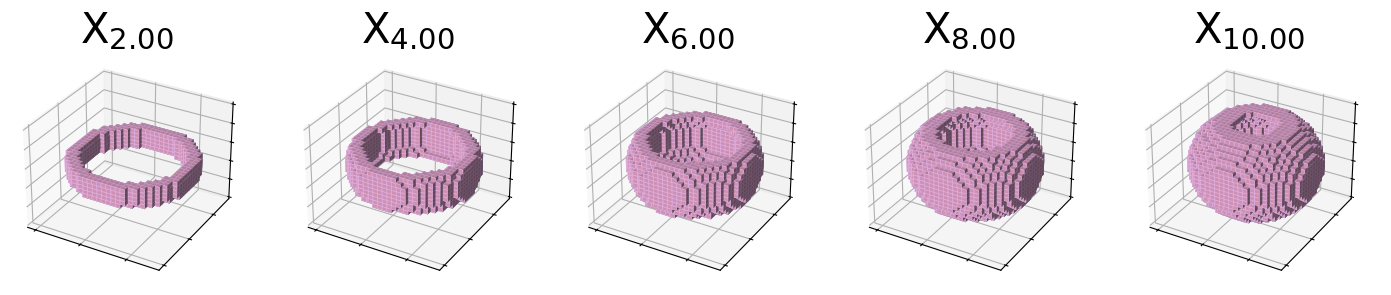} &
\includegraphics[height = 0.09\linewidth]{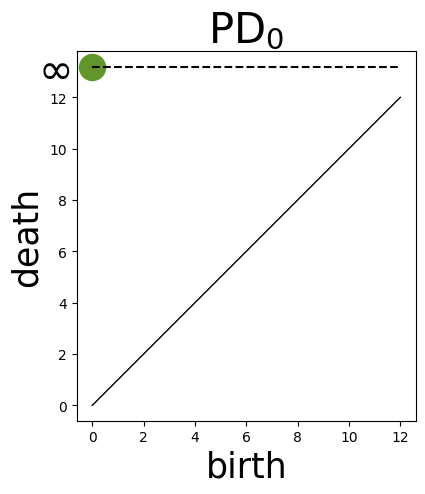} 
\includegraphics[height = 0.09\linewidth]{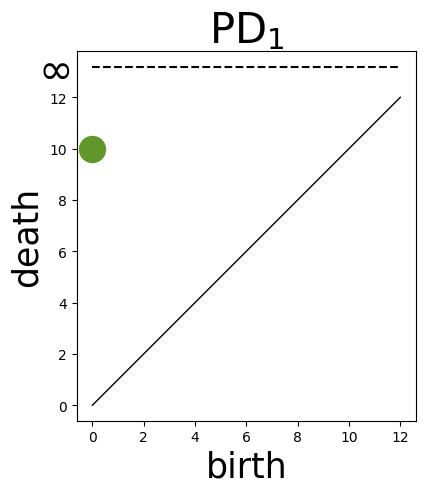}  
\includegraphics[height = 0.09\linewidth]{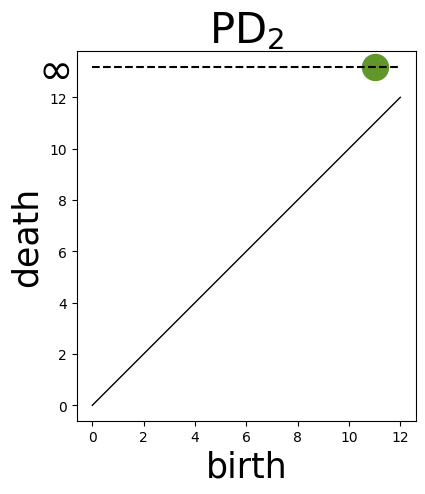} \\

\midrule

\includegraphics[height = 0.085\linewidth]{figures/tpht_vs_hpht/ball} &
\includegraphics[height = 0.085\linewidth]{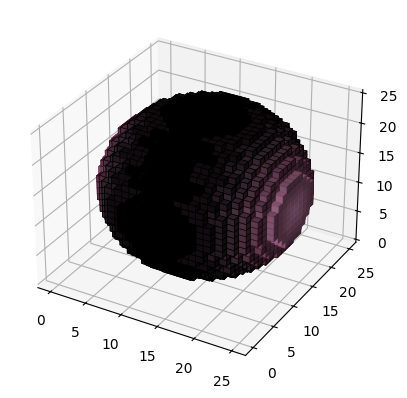} &
\includegraphics[height = 0.085\linewidth]{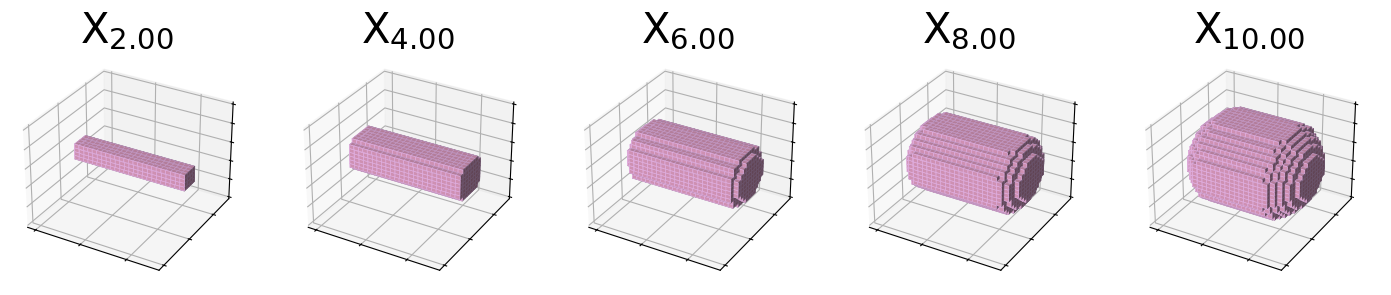} &
\includegraphics[height = 0.085\linewidth]{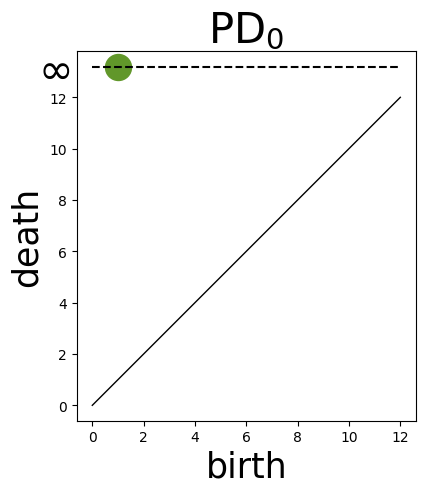} 
\includegraphics[height = 0.085\linewidth]{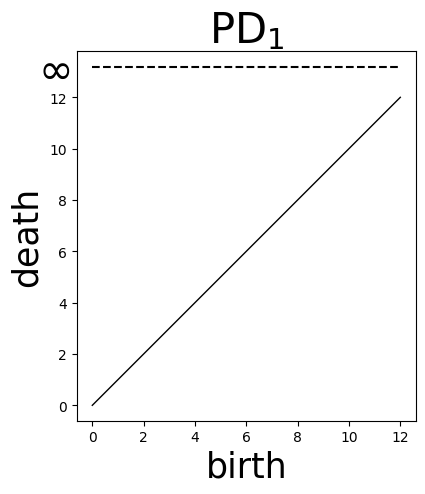} 
\includegraphics[height = 0.085\linewidth]{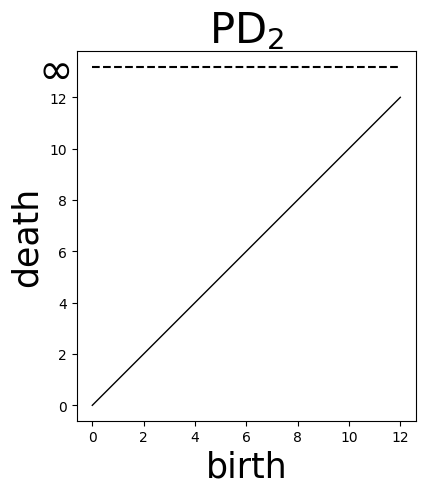} \\

\includegraphics[height = 0.085\linewidth]{figures/tpht_vs_hpht/sphere} &
\includegraphics[height = 0.085\linewidth]{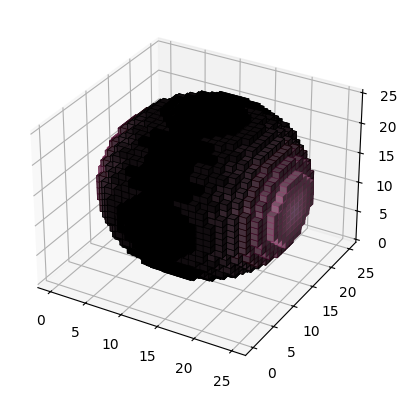} &
\includegraphics[height = 0.085\linewidth]{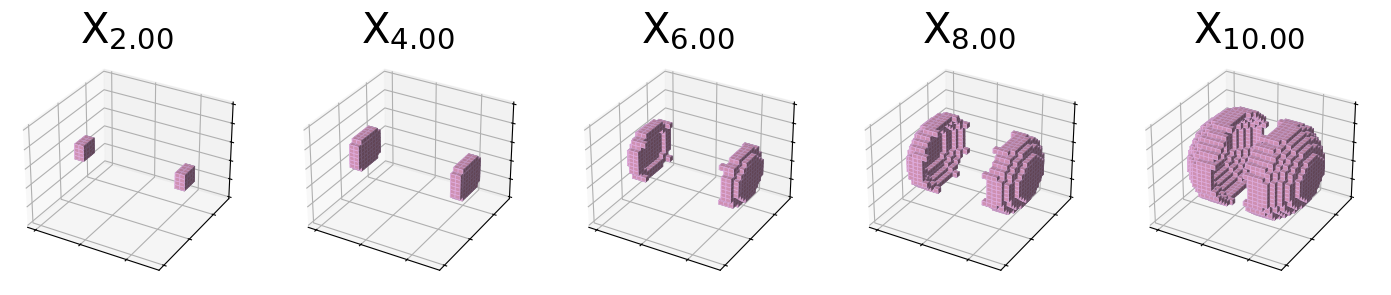} &
\includegraphics[height = 0.085\linewidth]{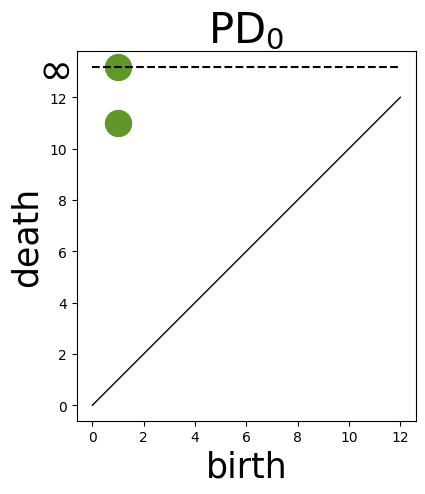} 
\includegraphics[height = 0.085\linewidth]{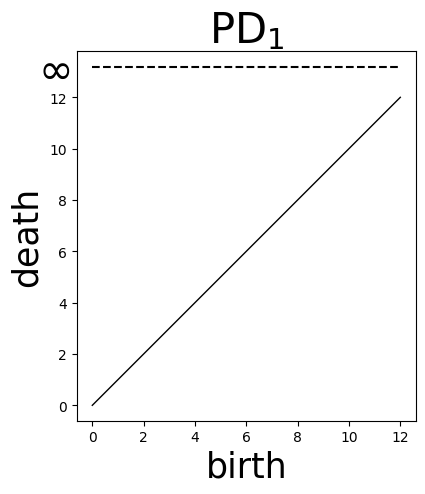}  
\includegraphics[height = 0.085\linewidth]{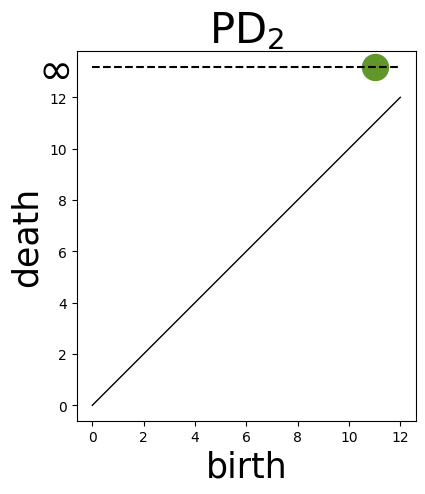} \\

\bottomrule
\end{tabular}

\caption{Illustration of an added value of distance-from-flat $\DPHT$ on $\AG(m', n)$ over $\AG(m, n)$ for $m' < m$, in particular of tubular $\TPHT$ over height $\HPHT$. $\PD_0(X)$ with respect to the adjusted height, i.e., distance from a plane (top rows) cannot discriminate between a ball- and sphere-like shape: there is always the one connected components in the filtration for both shapes; higher homological dimensions are needed. However, $\PD_0(X)$ with respect to the tubular filtration function, i.e., distance from a line (bottom rows) is sufficient to differentiate the two shapes: the sphere sees a second connected component in the filtration.}
\label{fig:tpht_vs_hpht}
\end{figure}

Our main contributions are as follows:

\begin{itemize}

\item We introduce generalised persistent homology transforms $\PHT_{\mathbb{P},f}$, allowing for arbitrary domain $\mathbb{P}$ and filtration functions $f_P,P\in \mathbb{P}$. This provides a broader framework for probing shapes, that opens doors to new research on $\PHT$ with respect to some other interesting families of filtration functions, that go beyond the classical height function.
\item We show that $\DPHT(X)$ is continuous (Theorem~\ref{thm:bottleneck-continuity}). In other words, if two flats (lines, planes, \dots, or hyperplanes) are close, then their corresponding persistence diagrams are close.

\item 
We show that $\DPHT$ \emph{truncated to homological degrees $k \in \{0, 1, \dots, m-1\}$}, is injective (Theorem~\ref{thm:injectivity}), which means that this information completely describes a shape (see Figure~\ref{fig:tpht_vs_hpht} for an illustration). 
Herein lies the main contribution of our proposed framework; in particular, $\TPHT$ in homology degree $0$ is a sufficient shape statistic --- this comes with significant computational advantage since it fairly easy to compute $\PH$ in degree $0$ in near-linear time with respect to the number $N$ of simplices by using union-find data structures \cite{edelsbrunner2010computational, wagner2012efficient, dlotko2014simplification}, whereas the standard algorithm is $\mathcal{O}(N^3)$ \cite{zomorodian2005computing}.

\item We show how  considering  the distance function $f_P=d_P$ with respect to {\em affine} subsets of $\mathbb{R}^n$ has an added value:
the classical $\PHT$ has a close relationship with $\DPHT$ for $m=n-1,$ since the \emph{absolute} value of the height filtration function $h_v$ can be seen as the distance from some hyperplane $P \in \AG(n-1, n)$; however, even in this case $\DPHT$ is more informative (see Proposition~\ref{proposition:pht_vs_hpht} and Figure \ref{fig:hpht_vs_cpht}).

\end{itemize}

To end, we note that the motivation for the generalised $\PHT,$ and $\DPHT$ in particular, came from an earlier work of some of the authors \cite{tubular}, where the tubular filtration is first introduced. There, we show that $\TPHT$ (that had not been formalised as a transform) in homological degree $0$ can detect convexity in $\mathbb{R}^n$ \cite[Theorem 1]{tubular}, and demonstrate its effectiveness and efficiency on synthetic and real-world data, where it outperforms some state-of-the-art deep learning techniques. We point interested readers to \cite[Remark 2]{tubular}, and in particular Figure 9, that clearly illustrates the added value of distances from affine hyperplanes compared to height, and distance from a line compared to a distance from a hyperplane (our two crucial generalisation aspects, or the last two items above), and has therefore served as the main catalyst for this work.
The remainder of the paper starts with the background (Section~\ref{sec:background}) needed to define our generalised persistent homology transform (Section~\ref{sec:gpht}), and to show its continuity (Section~\ref{sec:continuity}) and injectivity (Section~\ref{sec:injectivity}). Section~\ref{sec:conclusions} summarises the conclusions and some directions for further work.

\section{Preliminaries}
\label{sec:background}

In this section, we provide the background needed to define our generalised persistent homology transform $\GPHT$ (Section~\ref{sec:gpht}), and to prove our main results, namely the continuity of $\DPHT(X)$ (Section~\ref{sec:continuity}), and the injectivity of $\DPHT$ (Section~\ref{sec:injectivity}) on the affine Grassmannian.
We start with the preliminaries on the domain of $\DPHT(X)$ --- the affine Grassmannian (Section~\ref{sec:background-ag}) ---, and its codomain --- persistent homology (Section~\ref{sec:background-ph}).
The injectivity of $\DPHT$ will follow from the injectivity of a related transform based on Euler integrals (a so-called ``Radon transform''), 
which we show using Schapira's inversion formula.  
We therefore end this section with some background on o-minimal structures, providing the right finiteness properties of topological spaces for Euler calculus (Section \ref{sec:o-minimal}), the Euler characteristic and calculus (Section~\ref{sec:background-euler}) and Radon transforms (Section~\ref{sec:background-radon}).
For ease of reference, we provide a table summarising the notation used in Table \ref{tab:notation}.

\subsection{The Grassmannian of affine subspaces}
\label{sec:background-ag}

For $X \subset \mathbb{R}^n$, the domain of our persistent homology transform $\DPHT(X)$ of interest will be the affine Grassmannian $\AG(m, n),$ the set of affine $m$-dimensional subspaces in $\R^n$: the set of all points $(m=0)$,  lines $(m=1)$,  planes $(m=2)$, \dots, and  hyperplanes $(m=n-1)$.
Affine Grassmannians generalise the more commonly known Grassmannians, the set of all subspaces of a given dimension of Euclidean space, i.e., the set of all lines, planes, \dots, or hyperplanes that \emph{pass through the origin}. 
Affine Grassmannians are less well studied than Grassmannians, because of some inherent difficulties due to their non-compactness, or to the difficulty of defining a metric structure on them. Nonetheless, there has been recent advance in their study \cite{lim2021grassmannian}, which provides the basis for this section. In the following we recall basic facts about affine Grassmannians that will be needed to study our PHTs of interest. 
Since much of the properties of affine Grassmannians are derived by properties of Grassmannians, we recall basic properties of Grassmannians as well.

\subsubsection{A plurality of perspectives}

There are many different perspectives from which one can look at affine Grassmannians; in \cite{lim2021grassmannian} the authors provide a study of these spaces from the point of view of algebra, algebraic topology, algebraic geometry and differential geometry. 
The most useful perspective for our purposes is to model $\AG(m,n)$ as a vector bundle over $\Gr(m,n)$, and we discuss this briefly here. 
We first introduce some definitions and notation.

\begin{definition}[Grassmannian and affine Grassmannian spaces $\Gr(m, n)$ and $\AG(m, n)$] Let $n, m \in \mathbb{N}$ such that $m < n$. The \define{Grassmannian space} $\Gr(m, n)$ is the set of $m$-dimensional linear subspace of $\mathbb{R}^n$, referred to as \define{$m$-planes}. The \define{affine Grassmannian space} $\AG(m, n)$ is the set of $m$-dimensional affine subspace of $\mathbb{R}^n$, referred to as \define{$m$-flats}. 
\end{definition}


We will sometimes make use of the notation from \cite{lim2021grassmannian} to denote elements of affine Grassmannians. We usually denote an element of $\AG(m,n)$ by $P$, as we do for all elements of a parameter space of interest $\mathbb{P}$. When we need to specify relationships between the affine Grassmannians and the Grassmannian, we will write an element of $\AG(m,n)$ as     $P=\mathbb{A}+b$, where $\mathbb{A}$ is an $m$-plane in $\Gr(m,n)$ and $b\in \mathbb{R}^n$ is a vector in $\mathbb{R}^n$, which encodes the displacement of the $m$-flat $\mathbb{A}+b$ from the origin. Given an element $b=(b_1,\dots , b_n)^t$ of $\mathbb{R}^{n}$, we denote by $b'$ the element in $\mathbb{R}^{n+1}$, defined by $b'=(b_0,\dots , b_n, 0)^t$. Similarly, for a subspace  $\mathbb{A}$ of $\mathbb{R}^n$, and $A=\{a_1, \dots , a_m\}$ a basis of $\mathbb{A}$, we denote by $\mathbb{A}'$ the linear subspace of $\mathbb{R}^{n+1}$ spanned by the vectors $\{a_1', \dots , a_m'\}$.

There are two important maps relating affine Grassmannians to Grassmannians.

\begin{definition}
\label{D:maps aff gr}

The \define{deaffine map} 
\begin{alignat}{2}
\pi\colon \AG(m,n) &\notag\to \Gr(m,n)\\
\mathbb{A}+b &\notag\mapsto \mathbb{A}
\end{alignat}
sends an affine subspaces to the linear subspace corresponding to it. 

The \define{embedding}

\begin{alignat}{2}
j \colon \AG(m,n)&\notag\to \Gr(m+1,n+1)\\
\mathbb{A}+b&\notag \mapsto \mathrm{span}\left(\mathbb{A}'\cup \{b'+e_{n+1}\}\right)
\end{alignat}
where $e_{n+1}=(0,0,\dots , 0, 1)^t\in \mathbb{R}^{n+1}$ is the $(n+1)$th unit vector of $\mathbb{R}^{n+1}$.

\end{definition}

For example, the embedding map sends a point $(x,y)\in\R^2 \cong \AG(0,2)$ to the line through the origin spanned by $(x,y,1)$ in $\R^3$, or sends the line $y=r$ in $\AG(1,2)$ to the plane in $\R^3$ spanned by $(0,0,1)$ and $(0,r,1)$.
We make the following straightforward observation explicit for future reference:

\begin{lemma}
The deaffine map $p$ and embedding $j$ from Definition \ref{D:maps aff gr} are a surjection and inclusion of sets, respectively. 
\end{lemma}

\subsubsection{Distances and metrics on affine 
Grassmannians}\label{SSS:metr ag}
Metrics on affine Grassmannians are difficult to define.
While many of the perspectives with which one can define the affine Grassmannians admit various canonical metrics, this is often complicated by non-compactness and/or leads to somewhat trivial metrics.
See \cite[Section 8.1]{lim2021grassmannian} for a full discussion.

The solution is to define a metric on the Grassmannians, which is better studied, and then inherit a metric on affine Grassmannian through the embedding from Definition \ref{D:maps aff gr}.

\paragraph{Distances on Grassmannians}
Any distance on $\Gr(m,n)$ that is invariant with respect to rotations is defined in terms of principal angles \cite[Theorem 2]{YL16}.
Several of these distances appear in the literature and are summarised in \cite[Table 2]{YL16}.
Principal angles generalise the notion of an angle between lines to higher-dimensional planes which have angles of separation in multiple orthogonal directions.

Given two $m$-planes $\mathbb{A}$ and $\mathbb{B}$, let $M_\mathbb{A}$ and $M_\mathbb{B}$ be $m\times m$ matrices whose columns are orthonormal bases of $\mathbb{A}$ and $\mathbb{B}$ respectively.
We set $\theta_i(\mathbb{A},\mathbb{B})$ to be the $i^\mathrm{th}$ principal angle between $\mathbb{A}$ and $\mathbb{B}$ and $\sigma_i(M)$ to be the $i^\mathrm{th}$ singular value of a matrix $M$.
Explicitly, principal angles relate to singular values via the following equation
\begin{equation}
\cos\theta_i(\mathbb{A},\mathbb{B}) = \sigma_i(M_\mathbb{A}^T M_\mathbb{B}).
\label{eq:angle-singular}
\end{equation}

There is one distance on $\Gr(m,n)$, considered \emph{the} distance, and often also called \define{Grassmann distance}, which is the (unique) geodesic distance induced by the Riemmannian metric structure on $\Gr(m,n)$.
Let $\theta_i:=\theta_i(\mathbb{A},\mathbb{B})$ denote the $i^\mathrm{th}$ principal angle between $m$-planes $\mathbb{A}$ and $\mathbb{B}$.
The Grassmann distance $\dist_{\Gr(m,n)}$ is
\begin{equation}
\dist_{\Gr(m,n)}(\mathbb{A},\mathbb{B}) = \left(\sum_{i=1}^m\theta_i^2\right)^{1/2}
\label{eq:grassmann-distance}
\end{equation}

\paragraph{Distances on affine Grassmannians}
All Grassmannian metrics can be extended to affine Grassmannian metrics through the embedding from Definition \ref{D:maps aff gr}.

\begin{definition}
Fix a metric $\dist_{\Gr(m+1,n+1)}$ on $\Gr(m+1,n+1)$.
The metric $\dist_{\AG(m,n)}$ on $\AG(m,n)$ is the unique metric such that the embedding $j$ from Definition \ref{D:maps aff gr} is an isometry.
That is,
\[
\dist_{\AG(m,n)}(P,Q) := \dist_{\Gr(m+1,n+1)}(j(P),j(Q))
\]
\end{definition}

In particular, the Grassmann distance induces a distance on $\AG(m,n)$ which is the only possible distance on the affine Grassmannian compatible with the Riemmannina metric structure on $\AG(m,n)$ induced from that of $\Gr(k+1,n+1)$ (see Lemma 8.4 in \cite{lim2021grassmannian}).
We also note that Corollary 8.2 in \cite{lim2021grassmannian} provides an explicit formula for the geodesics connecting two points on $\AG(m,n)$.

\begin{remark}
Singular values are not only defined for square matrices, so principal angles between flats of different dimensions are also well-defined.
Thus, with some care, one can extend  $d_{\AG(m,n)}$ to measure distances between flats of different dimensions.
See Section 8.3 in \cite{lim2021grassmannian} for an explicit construction.
\end{remark}

\begin{remark}
The deaffine map $\pi$ is automatically continuous with respect to the vector bundle topology.
This is not the same topology as the metric space $(\AG(m,n),\dist_{\AG(m,n)})$, however in Section~\ref{sec:continuity} we will show that the deaffine map is also continuous with respect to this metric.
\end{remark}

\subsubsection{Schubert calculus}\label{SS: Schu}
Schubert calculus refers to enumerative techniques developed first by Hermann Schubert in his 1879 book \cite{schubert}, as a means to count the number of all possible subspaces  of a given dimension that intersect specific given subspaces.

\paragraph{The Grassmannian of all m-flats containing a specific k-flat}
A key ingredient of Schubert calculus are the Schubert varieties.
For a $k$-flat $\mathbb{A}+b$ in $\R^n$ with $0\leq j\leq k\leq l < n$, these are defined as
\begin{align*}
\Psi^+_l(\mathbb{A}+b) & := \left\{\mathbb{B}+c \in \AG(l,n)\,:\,\mathbb{A}+b\subseteq \mathbb{B}+c\right\} \\
\Psi^-_j(\mathbb{A}+b) & := \left\{\mathbb{C}+d \in \AG(j,n)\,:\,\mathbb{A}+b\supseteq \mathbb{C}+d\right\},
\end{align*}
and similarly we can define the Schubert varieties of planes when restricting to (non-affine) Grassmannians.
That is, the positive Schubert varieties of a flat (respectively, plane) are the all the higher-dimensional flats (planes) which contain it, and the other is all the lower-dimensional flats (planes) that it contains.
It is easy to see that Schubert varieties of a flat (plane) can in fact be identified with non-affine Grassmannians \cite{lim2021grassmannian} (\cite{YL16}).

\begin{proposition}
\label{prop:schubert}

If $\mathbb{A}+b\in\AG(k,n)$ then
\[
\Psi^+_l(\mathbb{A}+b)\cong \Gr(n-l,n-k)\;\;\;\;\;\;\text{and}\;\;\;\;\;\;\Psi^-_j(\mathbb{A}+b)\cong \Gr(j,k)
\]
\end{proposition}

One can use this to show, for instance, that the space of $k$-flats which connect two points in $\R^n$ can be identified with $\Gr(n-k,n-1)$.
In turn, identifying subspaces of $\R^n$ with their orthogonal compliment, we see that the space of $k$-flats connecting two points in $\R^n$ is homeomorphic to $\Gr(k-1,n-1)$.
We will explicitly use this fact in Section~\ref{sec:injectivity}.

\paragraph{Euler characteristic of Grassmannians} Schubert varieties are also a key ingredient in calculating the Euler characteristic of Grassmannians, which we require in Section~\ref{sec:injectivity}.
This construction is widely known, however we were unable to find any reputable references.
For completeness, we give a very brief construction of $\chi(\Gr(m,n))$ here and direct the reader to their favourite source of mathematical blogs for further details.

Let $\ell=\mathrm{span}\{e_1\}$ be a $1$-space in $\Gr(1,n)$.
Consider
\[
\Gr(k,n) = \Psi^+_k(\ell) \cup \Psi^+_k(\ell)^C.
\]
The Euler characteristic is additive, so observing that $\Psi^+_k(\ell)\cong \Gr(k-1,n-1)$ and that $\Psi^+_k(\ell)^C$ is a vector bundle over $\Gr(k,n-1)$ allows us to define a recursive relation
\[
\chi(\Gr(k,n)) = \chi(\Gr(k-1,n-1)) + (-1)^k\chi(\Gr(k,n-1))
\]
where $\Gr(0,n)\cong\R^n$ and $\Gr(n,n)\cong\{\star\}$ means $\chi(\Gr(0,n))=\chi(\Gr(n,n))=1$.
This has an explicit formula
\[
\chi(\Gr(k,n)) = \begin{cases}
0 & \text{ if }n\text{ even, and }k\text{ odd}\\
\binom{\lfloor n/2\rfloor}{\lfloor k/2 \rfloor} & \text{ otherwise}.
\end{cases}
\]
There are alternative constructions which also recover the Euler chracteristic, see for example \cite[Corollary 5.2]{casian2013cohomology}.

\subsection{Persistent homology}
\label{sec:background-ph}

In this section we introduce sublevel set persistent homology.
We note that while some  readers may be familiar with the persistent homology of distance-based filtrations of finite point clouds (e.g., persistent homology of filtrations of geometric complexes such Vietoris-Rips, \v Cech, alpha complexes), persistent homology can be defined for any filtration of a topological space.
For more details, we direct the reader to \cite{edelsbrunner2008persistent}.

\begin{definition}[Filtration]
\label{def:filtration}
A \define{filtration} $\mathcal{F}$ of a topological space $X$ is a nested family $\{\FF_r\}_{r \in \R}$ of subspaces of $X$. Of particular interest are sublevel sets of so-called \define{filtration functions} $f \colon X \rightarrow \R,$ where 
$$\FF_r = X_r^{-}(f) :=\{ x \in X \mid f(x) \leq r\}.$$
\end{definition}

Next to the sublevel sets, the related notion of level sets will be of interest, namely 
\[
X_r(f) =\{ x \in M \mid f(x) = r\}\, .
\]
As we discussed in the Introduction, in this paper we focus on $\PH$ on affine Grassmannians $\AG(m, n)$, with the distance $d_P$ from a flat $P \in \AG(m, n)$ as the filtration function, so in this case we simplify the notation:
$$X_r^- := X_r^-(d_P)  =\{ x \in X \mid d(x, P) \leq r\}$$
$$X_r := X_r(d_P) = \{ x \in X \mid d(x, P) = r\}.$$
Note further that $X_0 = \{ x \in X \mid d(x, P) = 0\} = X \cap P.$

Given a filtration $\mathcal{F}$ of a topological space $X$, we may compute its (singular) homology with coefficients in a field. We denote by $C_k(\FF_r)$ the vector space of singular $k$-chains and by $\partial_{k+1}\colon C_{k+1}(\FF_r)\to C_{k}(\FF_r)$ the  boundary operator defined by sending each singular simplex to the alternating sum of its faces.
Then, given 
  a filtration $\{\FF_r\}_{r \in \R}$,
we obtain injective linear maps $C_k(\FF_r)\hookrightarrow C_k(\FF_{r'})$ induced by the inclusion maps $ \FF_r\hookrightarrow \FF_{r'}$.  
Therefore, we may identify each $C_k(\FF_r)$ with a vector subspace of $C_k(\FF_{r'})$ for all $r\leq r'$.

\begin{definition}
Given $c\in C_k(\FF_r)$ we define its \define{birth} to be

\[
\birth(c)=\min\left \{r\in \mathbb{R}\mid c\in C_k(\FF_r)\right \}
\]
and its \define{death}
\[
\death(c)=\min\left \{r\in \mathbb{R}\mid c\in \im(\partial_{k+1}(C_{k+1}(\FF_r))\right\}\, ,
\]
with the convention $\min\emptyset =\infty$.

\end{definition}

 Similarly as for the singular chain complexes, the inclusion maps $\iota_{r,r'}\colon \FF_r\hookrightarrow \FF_{r'}$ induce maps on homology vector spaces $H_k(\iota_{r,r'})\colon H_k(\FF_r)\to H_k(\FF_{r'})$ satifying the conditions $H_k(\iota_{k,k})=id$ and $H_k(\iota_{r',r''})\circ H_k(\iota_{r,r'})=H_k(\iota_{r,r''})$ for all $r\leq r'\leq r''$. One defines:

\begin{definition}
The \define{$k$th persistence homology module} of a space $X$ with respect to a filtration  $\mathcal{F}=\{\FF_r\}_{r\in \mathbb{R}}$ is the tuple $\Big (\{H_k(\FF_r)\}_{r\in \mathbb{R}}, \{H_k(\iota_{r,r'})\}_{r\leq r'}\Big)$, denoted in short by $H_k(X,\mathcal{F})$.
\end{definition}

It is a fundamental theorem in Topological Data Analysis that for (singular) homology with coefficients in a field, and for $X$ satisfying appropriate finiteness conditions, there exists a choice of basis vectors of $C_k(X)$, called ``persistent homology cycle basis'',  which is compatible with the filtration in a precise sense. We refer the reader to \cite{LO24} for details on the definition of such a basis. Given such a persistent homology cycle basis, we can define the persistence diagram of $H_k(X,\mathcal{F})$ as the multiset of birth-death values for each vector in the basis:

\begin{definition}\label{def:pd}

Given a filtration $\mathcal{F}=\{\FF_r\}_{r\in \mathbb{R}}$ of a topological space $X$, and 
given   a persistent homology cycle basis $c_1,\dots , c_m\in C_k(X)$  of $H_k(X,\mathcal{F})$, we call the multiset
\[
PD_k(X,\mathcal{F}):=\Big\{\big(\birth(c_j),\death(c_j)\big) \mid  j=1,\dots , m \Big\},\]
together with the set of points on the diagonal with infinite multiplicity, the persistence diagram of $H_k(X,\FF)$, or the \define{persistence diagram of $X$ in homology degree $k$} (with respect to the filtration $\FF$).
We call the number $\death(c_j)-\birth(c_j)$ the \define{persistence} of $c_j$. 
If $\mathcal{F}$ is the sublevel set of a filtration function $f$, we
write $\PD_k(X, \mathcal{F}) = \PD_k(X, f)$.

Finally, we say that $c_i$ is an \define{essential cycle} if it has infinite persistence, and similarly, we call the pair $(\birth(c_i),\death(c_i))$ of such a cycle an \define{essential point} in the persistence diagram.
\end{definition}

For an illustration of a filtration and the resulting persistence diagrams, see, for example, the first row in Figure~\ref{fig:pht}: the essential persistence points in $\PD$ in homology degrees $0$ and $1$ capture respectively the one connected component and one hole, and the birth values capture their position.

\begin{remark}
For the persistence diagram of the sub-level sets of $f$ to be well-defined, one usually requires a mild regularity condition that $f$ is \define{tame}.
We say $f$ is tame if $H_k(X_r^-)$ is finite-dimensional for all $k\in\mathbb{Z}$ and $r\in\R$ and there are only finitely many homological critical values where $H_k(X_{r-\epsilon}^-)\to H_k(X_{r+\epsilon}^-)$ is not an isomorphism for sufficiently small $\epsilon$ \cite{CSEH06}.
This is slightly different from the notion of tameness in Section~\ref{sec:o-minimal}, although sufficiently nice functions on o-minimal sets will be tame.
\end{remark}

One of the selling points of persistent homology is that it is known to be stable under perturbations in many settings.
In order to quantify this, we need a notion of distance between persistence diagrams.

\begin{definition}[Wasserstein and bottleneck distance between persistence diagrams \cite{skraba2020wasserstein}]
\label{def:wasserstein_distance}

Let $\PD, \PD'$ be two persistence diagrams, and let $p , q \in [1, +\infty).$ The $p$-th \define{Wasserstein distance between the persistence diagrams} is defined as:
$$
W_{p, q}(\PD, \PD') = \inf_\tau \Big( \sum_i \| (b_i, d_i) - \tau(b_i, d_i) \|_q^p \Big)^\frac{1}{p},
$$
where the infimum is taken across all bijections $\tau: \PD \rightarrow \PD',$ and the sum across all persistence intervals $(b_i, d_i) \in \PD$. 
The \define{bottleneck distance} for $p=+\infty$ is defined as
$$
W_{\infty, q}(PD_1, PD_2) = \inf_\tau \sup_i \| (b_i, d_i) - \tau(b_i, d_i) \|_q.
$$
One commonly assumes that $q = \infty,$ and denotes $W_p = W_{p, \infty}$ for any $p \in [1, +\infty]$ and $d_B=W_{\infty,\infty}.$ Note that there exists a bijection between any two $\PD$s, since a $\PD$ includes the diagonal with infinite multiplicity (Definition~\ref{def:pd}).
\end{definition}

\begin{theorem}[Bottleneck stability \cite{CSEH06}] 
\label{thm:bottleneck-stability}
Let $X$ be triangulisable, and $f, g:X \to \R$ continuous tame functions. The bottleneck distance between $\PD_k(X,f)$ and $\PD_k(X,g)$ is bounded by the supremum distance between $f$ and $g$:
$$
d_B(\PD_k(X,f),\PD_k(X,g)) \leq \norm{f-g}_\infty.
$$
\end{theorem}

\begin{theorem}[Wasserstein stability \cite{cohen2010lipschitz}] 
\label{thm:wasserstein-stability}

Let $X$ be triangulisable, compact metric space that implies bounded degree-$p'$ total persistence: $\sum \{ (d_i-b_i)^{p'} \mid (b_i, d_i) \in \PD_k(X, \alpha)\} \leq C_X$ for every tame function $\alpha: X \to \R$ with Lipschitz constant  $\mathrm{Lip}(\alpha) \leq 1$. Let $f, g: X \to \R$ be tame Lipschitz functions. The $p$-Wasserstein distance between the corresponding persistence diagrams is bounded by 
the supremum distance between $f$ and $g$ up to a constant:
$$
W_p(\PD_k(X,f),\PD_k(X,g)) \leq C^{1/p} \norm{f-g}_\infty^{1-\frac{p'}{p}}.
$$
for every $p \geq k$, where $C = C_X \max \{ \mathrm{Lip}(f)^{p'}, \mathrm{Lip}(g)^{p'} \}.$

\end{theorem}

\subsection{A tameness framework for topology: o-minimal structures}\label{sec:o-minimal}

In his \emph{Sketch of a program} \cite{esquisse}, Grothendieck discusses the need for a ``tame topology'' (``topologie modérée''): the definition of a topological space is so general that it allows for wild behaviour that makes one's life unnecessarily complicated, if  one's aim is to study properties of geometric shapes invariant up to homeomorphism or homotopy.
Such a tame topology is provided by o-minimal structures, of which the ideas are already present in Grothendieck's program, and which have been formalised by model theorists in the 1980s, to generalise properties of semi-algebraic sets of $\mathbb{R}^n$, or subanalytic sets. 

Semi-algebraic sets are ``tame spaces''  because they are stable with respect to many desirable operations, such as taking unions, intersections, complements, projections or products. Thus, by taking repeated operations of this type one never encounters ``wild'' spaces, such as a topologist's sine curve.

Here we give a brief introduction to an o-minimal structure on the field of the reals $\mathbb{R}$, which provides  the right setting to compute integrals with respect to the Euler  characteristic. In turn, Euler calculus will be  a crucial ingredient in the investigation of the injectivity of our $\PHT$s. As references for o-minimal structures, we refer to the lecture notes in \cite{C20} and the monograph \cite{Dries_1998}.

\begin{definition}
A subset $A\subset \mathbb{R}^n$ is called \define{algebraic} if it can be defined by a finite number of polynomial equations with coefficients in $\mathbb{R}$:  given ${S=\{p_1,\dots , p_k\}\subset \mathbb{R}[x_1,\dots , x_n]}$, we can express $A$ as

\[A=\{(x_1,\dots , x_n)\in \mathbb{R}^n \mid p(x_1,\dots , x_n)=0 \;\; \forall p\in S \}\, .
\]

A subset $A\subset \mathbb{R}^n$ is called \define{semialgebraic} if it can be expressed as finite unions, intersections and complements of solution sets of polynomial equations or inequalities. 
\end{definition}

\begin{definition}\label{D:o-minimal}

An \define{o-minimal structure} on the field $\mathbb{R}$ is given by a collection of sets $\{S_n\}_{n\in \mathbb{N}_{\geq 1}}$ such that each $S_n$ is a subset of the powerset 
$\mathcal{P}(\mathbb{R}^n)$ and furthermore, these sets satisfy the following properties:
\begin{enumerate}[label=(\roman*{})]
\item each $S_n$ is closed under taking finite unions, intersections and complements.
\item if $A\in S_n$ and $B\in S_m$ then $A\times B\in S_{n+m}$.
\item if $A\in S_{n+1}$ and $\pi\colon \mathbb{R}^{n+1}\to \mathbb{R}^n$ is the projection onto the first $n$ coordinates, then $\pi(A)\in S_n$.
\item All algebraic subsets of $\mathbb{R}^n$ are elements of $S_n$.
\item The elements of $S_1$ are  the finite unions of points and intervals.
\end{enumerate}
 
 The elements of $S_n$ are called \define{definable} subsets of $\mathbb{R}^n$. A definable subset $X\subset \mathbb{R}^n$ is \define{constructible} if it is compact.  
We denote the set of constructible subsets of $\mathbb{R}^n$ by $CS(\mathbb{R}^n)$.

\end{definition}

\begin{example}
[Semialgebraic subsets  are definable] Consider the semialgebraic subset $S$ defined by the inequality $p(x_1,\dots , x_n)>0$. Then we have that $S=\pi (A)$ where $A\subset \mathbb{R}^{n+1}$ is the algebraic set defined by $x_{n+1}^2p(x_1,\dots , x_n)-1=0$ and $\pi$ is the projection onto the first $n$ components. Thus, $S$ is definable, and therefore the same is true of finite unions, intersection, complements of  subsets defined by finitely many polynomial equations and inequalities. 
\end{example}

\begin{example}[The affine Grassmannian is definable]\label{E: aff grass def}

Given $P=\mathbb{A}+b\in \AG(m,n)$, its \define{projection affine coordinates} are given by $(A,b_0)\in \mathbb{R}^{n\times (n+1)}$, where $A$ is a symmetric idempotent matrix with trace equal to $k$ (or, equivalently, of rank equal to $k$), and
 $b_0\in \mathbb{R}^n$ is a vector orthogonal to $\mathbb{A}$ and such that $\mathbb{A}+b=\mathbb{A}+b_0$. Thus, there is  a bijection 
 \begin{equation}\label{E:ag def}
\AG(m,n)\cong  S:=\left\{(A,b)\in \mathbb{R}^{n\times (n+1)} \mid A^2=A=A^T, \; \; A^Tb=0, \;\; \mathrm{tr}(A)=k \right\},
 \end{equation}
see \cite[Proposition 5.2]{lim2021grassmannian} for more details. The set $S$ on the right-hand-side of Equation \refeq{E:ag def} is algebraic, since the conditions $A^2=A$, $A=A^T$, $A^Tb=0$ and $\mathrm{tr}(A)=k$ are all polynomials.
In particular, $\AG(m,n)$ is definable.

\end{example}

\begin{definition}

    Given $X\subset \mathbb{R}^n$, we say that a function $f\colon X\to \mathbb{R}^m$ is \define{definable} if its graph is a definable subset of $\mathbb{R}^{n}\times \mathbb{R}^{m}$.

\end{definition}

\begin{example}[The distance-from-flat function is definable]\label{E: distance def}
Let $X\subset \mathbb{R}^n$ be constructible, $P\in \AG(m,n)$, and consider  the function
\[\dist(x, P):=\min_{p\in P} \| x - p \|_2 \, ,
\]
encoding the distance from $m$-flats $P \in \AG(m,n)$.
 Given $r\in \mathbb{R}$, we consider the sublevel set $X^{-}_r(\dist_P)$ and express it as the intersection of $X$ together with the sublevel sets ${\mathbb{R}^n}^{-}_r(\dist_P)$ of the function $d_P\colon \mathbb{R}^n\to \mathbb{R}$ sending each point in $\mathbb{R}^n$ to its distance from the flat $P$.

First, let $m=0$. Then  ${\mathbb{R}^n}^{-}_r(\dist_P)$ is an $n$-dimensional closed ball of radius $r$ centered at $P$. Such a ball can be expressed as the solution set of the polynomial inequality  $\sum_{i=1}^n(x_i-P_i)^2\leq r^2$, and therefore is semialgebraic, compact and in particular constructible. Therefore the intersection $X^{-}_r(\dist_P)=X\cap  {\mathbb{R}^n}^{-}_r(\dist_P)$ is constructible. 

 Now let $m=1$. For simplicity, assume first that $P$ is the line spanned by the vector $e_1$. Then ${\mathbb{R}^n}^{-}_r(\dist_P)$ is given by all points $(x_1,\dots , x_n)$ within distance $r$ from some point on the line. 
 Thus, the sublevel set consists of all points $(x_1,\dots , x_n)\in \mathbb{R}^n$ such that $x_2^2+\dots + x_n^2\leq r^2$. Hence ${\mathbb{R}^n}^{-}_r(\dist_P)$  can be identified with the product $P\times B$ where $B=\{(0,x_2,\dots , x_n)\mid x_2^2+\dots + x_n^2\leq r^2\}$ is a ball of dimension $n-1$. Thus, since $P$ is algebraic and $B$ is semialgebraic, $P\times B$ is semialgebraic, and thus in particular definable. Therefore the intersection ${\mathbb{R}^n}^{-}_r(\dist_P)=X\cap (P\times B)$ is definable. A similar argument shows that for arbitrary line $P\in \AG(1,n)$ the sublevel set ${\mathbb{R}^n}^{-}_r(\dist_P)$ can be expressed as $X\cap (P\times B)$ where $B$ is an $(n-1)$-dimensional closed ball.
 
A similar argument allows to deduce the following characterisation for arbitrary $m<n$. We have that 
\[
X^{-}_r(\dist_P)\cong X\cap  B_r
\]
where $B_r$ is a closed $(n-m)$-dimensional  ball of radius $r$.

In particular, the graph of $d_P$ can be written as
\[
\bigcup_{r\in \mathbb{R}}
\{
(x,r)\mid x\in X\cap S_r
\}
\]
where $S_r$ is a sphere of radius $r$ and dimension $n-m$. Since $X$ is compact, there exists $r_M$ such that $X\subset B_{r_M}$. Thus, we can rewrite the graph as 
\[
\bigcup_{r\in [0,r_M]}
\{
(x,r)\mid x\in X\cap S_r
\}
\]
Furthermore, since the homeomorphism type of the (sub)level sets changes only finitely many times (see Lemma \ref{L:finitely many changes}), one may identify the graph with the finite union 
\[
(X\cap S_{r_0'})\times [0,r_0]\cup \cdots  \cup  (X\cap S_{r_k'})\times [r_k,r_M]\]
for $r_0<\dots < r_k\in [0,r_M]$. 
Thus $d_P$ is definable.
We note that the description of the sublevel sets given here is suggestive of another interpretation of such sublevel sets as  ``discal'' filtrations. 
\end{example}

One of the main tools in the study of definable sets is the fact that every definable set has a finite partition into definable subsets $C_i$, each homeomorphic to an open box (i.e., a cartesian product of open intervals), see \cite[Theorem 2.10]{C20}. Such a partition is called ``cylindrical definable cell decomposition'', and the sets $C_i$ are called ``cells''.  For instance, a cylindral definable cell decomposition of $\mathbb{R}$ is given by a choice of finitely many points $a_1<\dots < a_k$. The cells are the singletons $\{a_i\}$ together with the open intervals delimited by these points.

If additionally we assume that the definable subset is compact, we have the more common description via a triangulation: any constructible set has a finite triangulation. More precisely, one has the following:

\begin{theorem}{\cite[Theorem 4.4]{C20}, \cite[Chapter 8]{Dries_1998}}
Let $X\subset \mathbb{R}^n$ be constructible. Let $C_1, \dots , C_k$ be definable subsets of $X$. 
There exists a finite simplicial complex $K$ and a definable homeomorphism $\phi\colon |K|\to X$ such that each $C_i$ can be written as $\cup_{j=1}^\ell\phi(\sigma_j)$ where $\sigma_j$ are open simplices of $K$.
\end{theorem}

\subsection{Euler characteristic and calculus}
\label{sec:background-euler}

The Euler characteristic is one of the simplest topological integer-valued invariants, which can be thought of as a generalisation of the cardinality of a finite set. For a topological space $X$ its \emph{homological} Euler characterisic $\chi_h(X)$ can be defined as the alternating sum of Betti numbers, whenever the sum is well-defined:
$$\chi_h(X) = \beta_0(X) - \beta_1(X) + \beta_2(X) - \beta_3(X) + \dots.$$ 

If one wishes to define an integration theory with respect to the Euler characteristic, one needs however to consider its combinatorial definition, which we give here in its generalisation to definable sets:

\begin{definition}
Let $X\subset \mathbb{R}^n$ be a definable set, together with a cylindrical definable cell decomposition $C_1,\dots , C_k$, where each cell $C_i$ is homeomorphic to $\mathbb{R}^k_i$. The \define{definable Euler characteristic} of $X$ is 
$\chi(X)=\sum_{i=1}^k(-1)^{\dim(C_i)}$.

\end{definition}

On constructible sets, the definable Euler characteristic agrees with the homological Euler characteristic.
However, unlike the homological Euler characteristic, the definable Euler characteristic is not a homotopy invariant of definable sets. For instance,  we have $\chi((0,1))=-1$ but $\chi([0,1])=1$. On the other hand, the definable Euler characteristic behaves like a measure, as $\chi(A\cup B)=\chi (A)+\chi(B)-\chi(A\cap B)$, for any definable subsets $A,B$.
As Propp discusses in \cite{propp}, there is a tension between the requirement that the Euler characteristic be homotopy invariant or that it behaves like a measure. 
The functions that we can integrate with respect to this measure are given by finite $\mathbb{Z}$-linear combinations of indicator functions:

\begin{definition}
Let $X$ be definable. A function $\phi\colon X\to \mathbb{Z}$ is called  \define{constructible} if every level set $\phi^{-1}(z)$ is definable, and only finitely many level sets are non-empty. We denote the set of constructible functions on $X$ by $CF(X)$.
\end{definition}

\begin{definition}[Euler integral \cite{curry2012euler}]
Let $U$ be a set, and $\phi\colon  U \rightarrow \mathbb{Z}$ a constructible  function. The \define{Euler integral  of $\phi$} is the sum of the Euler characteristic of each of the level sets of $\phi:$

$$\int_U \phi d\chi = \sum_{r = -\infty}^{r = \infty} r \cdot \chi(\{u \in U \mid \phi(u) = r\}).$$
\end{definition}

Taking the indicator function $\phi = \mathbf{1}_X$ for some $X \subset U$, it is easy to check that indeed $\int_U \mathbf{1}_X d\chi = \chi(X).$

\subsection{Radon integral transform}
\label{sec:background-radon}

The concepts discussed in this work are closely connected to those in integral geometry, which have been applied to model surfaces, random fields, point processes, but also in imaging, signal processing, and inverse problems \cite{ghrist2014elementary}. 
Integral geometry originated as an effort to refine aspects of geometric probability theory, building from the classic theorem of Crofton, which expresses the length of a plane curve as as an expectation of the number of intersections with a random line. 

The core idea in integral geometry is to examine invariant integral transforms, a specific type of mathematical operator that maps a function from its original function space into a different function space via integration
$$(T \phi)(v) = \int_{t_{1}}^{t_{2}} \phi(u) \, K(t,u) \, dt,$$
where some of the properties of the original function might be more easily characterised and manipulated than in the original function space (e.g., a function on a surface or shape is transformed into a function defined on an interval of the real line). There are many useful integral transforms, each specified by a choice of the function $K$ of two variables, referred to as the kernel or nucleus of the transform. Some well-known examples are the Fourier and Laplace transforms.

A perhaps less well-known integral transform is the Radon transform 
introduced by Radon in 1917 \cite{radon}, which sends a function defined on $\mathbb{R}^3$ to all of its line integrals. Later,  the Radon transform was generalised for arbitrary dimension $n$ \cite{rubin2004radon},  and for constructible subsets by Schapira \cite{schapira1995}. The Radon transform in $\mathbb{R}^3$ is commonly used in tomography, where images are generated from projection data obtained through cross-sectional scans of an object. 

Here, we will make use of the formulation in Schapira's work, since Schapira's inversion result will be crucial to prove the injectivity of $\DPHT$. 
In its general formulation, the Radon transform allows us to study a constructible function by sending it to the function that computes its Euler integral with respect to any constructible subset of $\mathbb{R}^n$.

\begin{definition}[Radon integral transform \cite{schapira1995}]
\label{def:radon}

Let $U$, $V$ be definable sets, and $S \subset U \times V$ a locally closed definable subset.
The \define{Radon transform} with kernel $S$ is the map $\radon_S : CF(U) \rightarrow CF(V)$ defined explicitly in the following way:

$$(\radon_S \phi)(v) = \int_{U} \phi(u) \, \mathbf{1}_S(u, v) \, d\chi(u).$$

\end{definition}

The principal result of Schapira \cite[Theorem 3.1]{schapira1995} gives a topological criterion for the invertibility of the Radon transform\footnote{We provide the formulation from [Theorem 2.13]\cite{CST22}}:

\begin{theorem}[Schapira's inversion formula \cite{schapira1995}]
\label{thm:schapira}

Let $U$, $V$ be definable sets and let $\phi \in \CF(U)$.
Let further $S \subset U \times V$ and $S' \subset V \times U$ be relations with fibers 
${S_u=\{u \in U: (u, v) \in S\}}$ and 
${S'_u=\{u \in U: (v, u) \in S'\}}$ 
satisfying the following conditions:
\begin{itemize}
\item[(1)] $\chi(S_u\cap S'_u) = \chi_1$ for all $u \in U$
\item[(2)] $\chi(S_u\cap S'_{u'}) = \chi_2$ for all $u \neq u'\in U$.
\end{itemize}
Then
$$(\radon_{S'} \circ \radon_{S})(\phi) = (\chi_1 - \chi_2) \phi + \chi_2 \left( \int_U \phi d\chi \right) \mathbf{1}_{U}.$$
\end{theorem}

In particular, note that this result implies that, when $\chi_1 \neq \chi_2,$ it is possible to recover $\phi$ from the inverse transform (and an affine map) \cite{GLH2019}. As Schapira concluded in his seminal paper \cite{schapira1995}, this inversion formula allows to reconstruct a body in a three-dimensional vector space from the integrals of the characteristic function of all its slices, that is, from the number of connected components minus the number of holes of all its intersections by two-dimensional affine slices.

\section{(Distance-from-flat) persistent homology transforms}
\label{sec:gpht}

We start this section by introducing the general persistent homology transform $(\PHT)$ that allows to study shapes through the lens of $\PH$ with respect to any filtration function. We then give the main definition of the paper, which formalises the distance-from-flat persistent homology transforms $\DPHT$, which constitute the focus of this work. In the remainder of the section, we discuss the relationship with the classical $\CPHT$ introduced in \cite{TMB2014} (Section~\ref{sec:gpht:subsec:cpht}), and some other examples that could be seen as special cases (Section~\ref{sec:gpht:subsec:examples}) or some extensions (Section~\ref{sec:gpht:subsec:other}).

\begin{definition}[$\PHT$]
\label{def:gpht}
Let $X \subset \R^n$ be a constructible set. Let further $\mathbb{P}$ be a topological space, and $\{f_P\}_{P \in \mathbb{P}}$ be a family of functions $f_P: X \rightarrow \mathbb{R}$. The \define{persistent homology transform ($\PHT$)} of $X$ with parameter space $\mathbb{P}$ and filtration functions $\{f_P\}_{P\in \mathbb{P}}$ is the function
\begin{alignat}{3}
\PHT(X) \colon & \notag \mathbb{P} && \notag \to \mathcal{D}^n\\ 
& P && \notag \mapsto (\PD_0(X, f_P), \PD_1(X, f_P), \dots, \PD_{n-1}(X, f_P) ) \, ,
\end{alignat}
where $\PD_k(X, f_P)$ is the persistence diagram of $X$ in homology degree $k$, with respect to the sublevel-set filtration of the function $f_P$ (Definition~\ref{def:pd}). When it is important to stress the choice of the parameter space and the filtrations, we write $\GPHT.$
\end{definition}

Two desirable properties of such $\PHT$s that may make them useful in shape analysis problems are continuity and injectivity: 

\begin{definition}[Injectivity of $\PHT$]
Given a $\PHT$ with parameter space $\mathbb{P}$ and filtration functions $\{f_P\}_{P\in \mathbb{P}}$, 
we say that $\PHT$ is \define{injective} if $\PHT(X) \ne \PHT(Y)$ for any constructible sets $X \ne Y$. 
\end{definition}
We further discuss continuity in Section~\ref{sec:continuity}, but note that it is different to say that the function $\PHT(X):\mathbb{P}\to\D^n$ is continuous than it is to say that the transform $\PHT$ is continuous.
The former has $\mathbb{P}$ as a domain, whereas the latter has the set of constructible functions as a domain, which will not necessarily be stable in terms of persistence.

\begin{remark}[Persistence bundles]
Our $\PHT$s are closely related to  persistence bundles introduced in \cite{hickok2022persistence}, however they differ from them in two main aspects. A persistence bundle is defined as a family of simplicial complexes $\{K_p\}_{p\in \mathbb{P}}$, together with a family of filtration functions $\{f_P\colon K_p\to \mathbb{R}\}$, together with the ``total space'', which summarises the information given by persistence diagrams in any homology degree for any filtration function. Of note is that persistence diagrams are assumed to contain only points of multiplicity $1$. Thus, on one hand persistence bundles are more general than $\PHT$s, since they also allow for a family of \emph{spaces} parametrised by the parameter space $\mathbb{P}$, while on the other hand our setting does not consider any restriction on the space of persistence diagrams. We note that while generically, if one computes persistent homology of point clouds, with respect to, for instance, a Vietoris-Rips or alpha complex, such an assumption is not too restrictive, the same is not true for other settings, including the one considered in the present  work, or  applications to networks.
\cite{hickok2022persistence}
\end{remark}

In this work we are interested in studying the following special cases of $\PHT$s:

\begin{definition}
The \define{distance-to-flat persistent homology transform} is the $\PHT$ where the domain  is the affine Grassmannian space, $\mathbb{P} = \AG(m,n),$ and the filtration functions $f_P(x) = \dist(x, P)$
encode the distance from $m$-flats $P \in \AG(m,n)$. We use $\DPHT$ to denote such transforms. 
As special cases of $\DPHT$ in this work we consider the distance from hyperplanes $(m=n-1)$, lines $(m=1)$ and points $(m=0)$, which we refer to  as the \define{height, tubular} and \define{radial persistent homology transform}, respectively. 
\end{definition}

\subsection{Relationship with the classical $\CPHT$}
\label{sec:gpht:subsec:cpht}

As we indicate in the Introduction, the classical $\CPHT$ introduced in \cite{TMB2014} is a $\PHT$ according to the above definition, for $\mathbb{P} = \Sp^{n-1}$ and $f_v(x) = x \cdot v$ the height function in direction $v$ (Figure~\ref{fig:pht}).
$\CPHT$ is famously injective, and known to be continuous with respect to the Bottleneck and Wasserstein distances on $\mathcal{D}^n$. While it is closely related to $\HPHT$, it is not a distance-from-flat PHT, as we formalise in the next proposition.

\begin{proposition}[Classical $\CPHT$ is less powerful than $\DPHT$.]
\label{proposition:pht_vs_hpht}
For any constructible $X \subset \R^n$, the image of the classical $\CPHT(X)$ can be identified with a subset of the image of $\HPHT(X)$. Furthermore, there exists some $X$ such that this subset relationship is strict. 
\end{proposition}

\begin{proof}

Let $X \subset \mathbb{R}^n$ be constructible, and let $B \subset \mathbb{R}^n$ be a ball of radius $M$ large enough so that $X$ is contained in its interior.
Let $X_r^-(h_v)$ be the sublevel set of the height function $h_v(x) = x \cdot v$ for some direction $v \in S^{n-1},$ and some $r \in \mathbb{R}.$ This set can be realised as the sublevel set of the distance function $d_P(x) = d(x, P),$ for appropriately chosen hyperplane $P \in \AG(n-1, n)$.
Indeed, for $P$ that is orthogonal to $v,$ and passes through the point in $B$ that is the furthest from $v,$ we have that $X_r^-(h_v)=X_{r+M}^-(d_P)$, since $d(x, P) = x \cdot v + M$, see Figure~\ref{fig:hpht_vs_cpht}(a).
Thus, we may identify the image of $\CPHT(X)$ with a subset of the image of $\HPHT(X)$ through the translation on the space of persistence diagrams $\mathcal{D}^n$, which takes an element $(\PD_0(X,h_v), \dots , (\PD_{n-1}(X,h_v))$ and translates every point in every diagram by $(M, M)$.
This completes the first part of the proof: $\HPHT$ is \emph{as informative} as $\CPHT.$

\begin{figure}[h!]
\centering
\begin{tabular}{ccc}
\includegraphics[height = 0.5\linewidth]{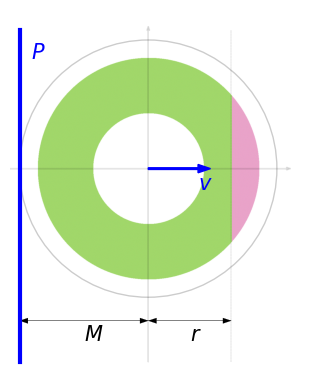} &&
\includegraphics[height = 0.485\linewidth]{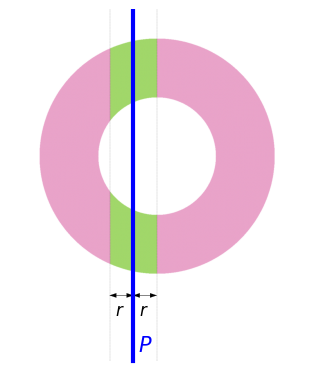} \\
(a) && (b) 
\end{tabular}
\caption{The height $\HPHT$ is more refined than the classical $\CPHT.$ Indeed, for any direction $v,$ the sublevel set (highlighted in green) of the height function $h_v$ can be realised as the sublevel set of the distance function $d_P$ for the appropriate hyperplane $P,$ $X_r^-(h_v) = X_{r+M}^-(d_P)$ (a). However, the sublevel sets of the distance from hyperplanes $P$ that pass through the shape (b) cannot be realised as sublevel sets of the height function $h_v$ for any direction $v$: the former can be seen as the region of the shape squished between two hyperplanes, whereas the latter is the intersection of the shape with a halfspace.} 
\label{fig:hpht_vs_cpht}
\end{figure}

To see why $\HPHT$ is \emph{more informative} than $\CPHT,$ note that the sublevel sets of the height and distance function are regions of the shape $X$ that lie within some halfspace, versus between two hyperplanes (Figure~\ref{fig:hpht_vs_cpht}). For a shape $X$ such as the annulus, the sublevel set of the distance from a hyperplane $P$ that passes through the hole in the shape consists of two connected components, highlighted in green in Figure~\ref{fig:hpht_vs_cpht}(b). This means that there are two points in $\PD_0(X, d_P),$ whereas there is a single point in $\PD_0(X, h_v)$ for any direction $v$. Therefore, the image of $\CPHT(X)$ is a strict subset of the image of $\HPHT(X).$

\end{proof}

In other words, the classical $\CPHT$ can be seen as the height $\HPHT$ restricted to the tangential hyperplanes in the above proof. And indeed, as we show later (Theorem~\ref{thm:injectivity}), the height $\HPHT$ can completely describe the shape with $\PH$ in homological degrees $0, 1, \dots, n-2,$ whereas the classical $\CPHT$ requires all degrees $0, 1, \dots, n-1.$

\subsection{Further examples of $\GPHT$}
\label{sec:gpht:subsec:examples}

\begin{example}[$\GPHT$ on hypersurfaces] 
\label{ex:gpht:hypersurface}
Although we focus on $\GPHT$ with $\mathbb{P}=\AG(m, n)$ --- the set of all flats (lines, planes, \dots, and hyperplanes), it is also possible to consider a more general case of hypersurfaces,
see \cite{ji2023}, in which the authors introduce an Euler characteristic transform with respect to quadric hypersurfaces, namely a hypersurface that can be written as $x^t Ax+v\cdot x=t$ where $A$ is a symmetric $n\times n$ matrix and $x\in \mathbb{R}^n$. Thus, we could define a $\GPHT$ with $\mathbb{P}$ given by the set of quadric hypersurfaces, and $f_P$ the distance from such a hypersurface. We note that the authors of \cite{ji2023} show that this transform is injective, and motivate considering such hypersurfaces by the detection of curvature, however, they don't explore this possibility further. 

\end{example}

\begin{example}[Filtrations on images] 
\label{ex:gpht:images-height}
In \cite{K19} persistent homology and tropical coordinates on the space of persistence barcodes are used to classify  digits in the MNIST data set. First, each grey-scale image is converted into a binary image, and four different filtrations are built on a cubical complex associated to the image. These four filtrations correspond to swiping the image from left-to-right or right-to-left, or from top-to-bottom and bottom-to-top. Considering four filtrations allows to retain spatial location of topological features in the images. We can thus consider this as an example of a $\PHT$, in which $\mathbb{P}=\{(1,0), (0,1), (0,1), (0,-1)\}$, and $f_p=h_p$. In particular, we would expect that by increasing the number of directions considered one would obtain a better classification accuracy. 
\end{example}

\begin{example}[$\GPHT$ on an arbitrary family of filtrations]
\label{ex:gpht:images-height-greyscale}
In the above examples, the family of filtrations $\{f_P\}_{P \in \mathbb{R}}$ are related so that each $f_P$ is defined in the same way, depending on some underlying parameter (direction $v$ for the height filtration, or flat $P$ for the distance-to-flat $\PHT$). We note that this does not need to be the case, since the definition of $\GPHT$ is very general.

For instance, in the example of digital images from Example \ref{ex:gpht:images-height}  one could take 
$\mathbb{P}=\{1,2,3,4,5\}$ where $f_1,\dots , f_4$ are the four filtrations described in Example \ref{ex:gpht:images-height}, and $f_5\colon C\to \{0,\dots, 256\}$ is the sublevel set filtration on the cubical complex  $C$ associated to the image that sends each vertex to the gray-scale value of the pixel it represents, edges between vertices to the maximum of the grey-scale values of the adjacent pixels, and squares to the maximum of the values of their boundary edges. 
\end{example}

\begin{example}[Persistence vineyards]
In the setting of persistence vineyards \cite{CSEM06}, one considers $X\subset \mathbb{R}^n$ triangulisable and a given abstract simplicial complex $K$ such that 
$|K| \cong X$, where $K$ is a geometric realisation of $K$, together with a $1$-parameter family of filtrations on $K$, encoded by a function $f_t\colon K\to \mathbb{R}$ for $t\in I\subset \mathbb{R}$. Thus, a vineyard is a $\PHT$ with $\mathbb{P}=I$ and $f_p=f_t$.
\end{example}

\subsection{$\PHT$ on non-Euclidean shapes}
\label{sec:gpht:subsec:other}

 We note that in Definition~\ref{def:gpht} we  require that a PHT is defined for a constructible subset of $\mathbb{R}^n$. Thus, considering different types of o-minimal structures or tameness from the one considered here, one may be able to consider $\PHT$s defined for sets $X$ that are not subsets of Euclidean space.

Some distance-to-flat PHTs naturally generalise to these non-Euclidean settings.
For instance, if $X\subseteq Y$ then the radial PHT of $X$ will be defined whenever $Y$ embeds \emph{definably} into Euclidean space.
Otherwise, if $Y$ is a geodesic metric space, then the tubular PHT of $X$ can be defined by replacing $\AG(1,n)$ with the space of all geodesic curves and filtering by the minimum distance to the geodesic. As we have shown in Example \ref{E: distance def} for our distance-from-flat function, a key ingredient in being able to define $\PHT$s is the definability of the filtration functions one considers.
Here we give some examples of potential candidates for non-Euclidean $\PHT$s for possible exploration, and we leave the question of whether these are true $\PHT$s (i.e., the persistence diagrams exist for all sublevel sets) for future work\footnote{In absence of a well-defined persistence diagram, one may   consider transforms mapping into the space of persistence modules instead of persistence diagrams.}.

\begin{example}[Tubular filtrations on the sphere]
Consider the $n$-sphere $\Sp^n$ and let $\mathcal{G}(\Sp^n)$ denote the geodesics of $\Sp^n$.
$\Sp^n$ is a definable subset of $\R^{n+1}$, so $\Sp^n$ inherits the Euclidean o-minimal subsets.
Explicitly, $\mathcal{G}(\Sp^n)$ contains the great circles of $\Sp^n$ (the isometric copies of $\Sp^1$).
For $X\subseteq \Sp^n$, the tubular/distance-to-geodesic PHT, $\PHT_{\mathcal{G}(\Sp^n),\dist}$, is the function sending $G\in \mathcal{G}(\Sp^n)$ to the persistence diagram with respect to the sublevel sets of $f_G:X\to \R$, the geodesic distance from $G$.

Consider the antipodal map $a:\Sp^n\to\Sp^n$ defined so that $a(x)=-x$ and observe that, since geodesics of $\Sp^n$ are great arcs, this means $\PHT_{\mathcal{G}(\Sp^n),\dist}(X)=\PHT_{\mathcal{G}(\Sp^n),\dist}(a(X))$
In particular, $\PHT_{\mathcal{G}(\Sp^n),\dist}$ is not injective.
\end{example}

\begin{example}[Tubular filtrations on the torus]
Consider the $n$-torus $\mathbb{T}^n$, let $\mathcal{G}(\mathbb{T}^n)$ denote the geodesics of $\mathbb{T}^n$, and fix some o-minimal structure on $\mathbb{T}^n$.
Modeling $\mathbb{T}^n=(\R/\mathbb{Z})^n$, $\mathcal{G}(\mathbb{T}^n)$ will be the straight lines of $\R^n$ (i.e. $\AG(1,n)$) modulo 1 in each coordinate.
For $X\subseteq \mathbb{T}^n$, the tubular/distance-to-geodesic PHT, $\PHT_{\mathcal{G}(\mathbb{T}^n),\dist}$, is the function sending $G\in \mathcal{G}(\mathbb{T}^n)$ to the persistance diagram with respect to the sublevel sets of $f_G:X\to \R$, the geodesic distance from $G$.

In this case, observe that lines in $\mathbb{T}^n$ of irrational slope may be dense in $\mathbb{T}^n$ but lines of rational slope will never be dense. Thus, it might be the case that the persistence diagram will not be defined for some sublevel sets, and in such case one might instead be interested in studying more generally ``persistence modules transforms'', sending a constructible set to the $n$-fold product of the space of persistence modules. 
The fact that rationals and irrationals are dense in $\R$ means that for $X\subseteq \mathbb{T}^n$ it is possible that $\PHT_{\mathcal{G}(\mathbb{T}^n),\dist}(X)$ is nowhere-continuous.
\end{example}

\begin{example}[Tubular filtration on the hyperbolic space] \label{ex:hyperbolic}
Consider the $n$-dimensional hyperbolic space $\Hp^n$ and let $\mathcal{G}(\Hp^n)$ denote the geodesics of $\Hp^n$.
Modeling $\Hp^n$ as the upper half-plane $\R^{n}\times (0,\infty)$, $\Hp^n$ inherits the o-minimal structure of $\R^{n+1}$ and $\mathcal{G}(\Hp^n)$ will be the union of the collection of circles centered on the hyperplane $x_{n+1}=0$ and the collection of vertical lines in $\R^{n+1}$.
For $X\subseteq \Hp^n$, the tubular/distance-to-geodesic PHT, $\PHT_{\mathcal{G}(\Hp^n),\dist}$, is the function sending $G\in \mathcal{G}(\Hp^n)$ to the persistance diagram with respect to the sublevel sets of $f_G:X\to \R$, the geodesic distance from $G$.
\end{example}

In Section~\ref{sec:injectivity} we will show that $\PHT_{\AG(m,n),\dist}$ is an injective statistic;  a similar argument shows that $\PHT_{\mathcal{G}(\Hp^n),\dist}$ is also a sufficient statistic.
We also note that it is not clear how, if at all, the classical $\PHT$ extends to hyperbolic space or other manifolds.
This shows that even aside from computational considerations, one benefit of our framework is that we can now introduce an injective statistic of shapes in hyperbolic spaces, and provide a framework for injective shape statistics on other manifolds where the tubular $\PHT$ may not be injective.

\begin{example}[Weighted networks]
\label{E:networks}
In \cite{KDS18} the authors study trees arising from brain networks by computing persistent homology of radial filtrations, from a designated root of the tree. More generally, given a weighted network, one might be interested in studying radial filtrations from a given finite subset of vertices in the network. This might for instance help in quantifying self-similarity of networks, by comparing the rate of growth of persistence intervals in persistence diagrams associated to different vertices, as one rescales the weights of the  network.
\end{example}

\begin{figure}[h!]
\centering
\scalebox{0.7}{
\begin{tabular}{c|c|c|c}
\toprule
Shape $X$ & $f_P$ & Filtration & $\PHT(X)$ \\
\midrule

\includegraphics[height = 0.15\linewidth]{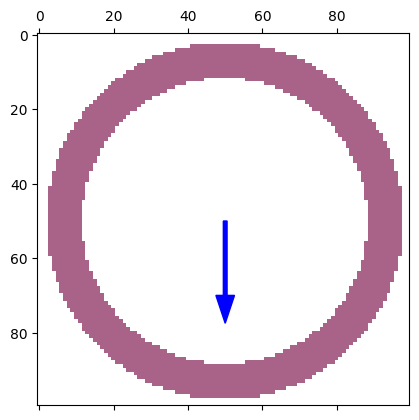} & 
\includegraphics[height = 0.15\linewidth]{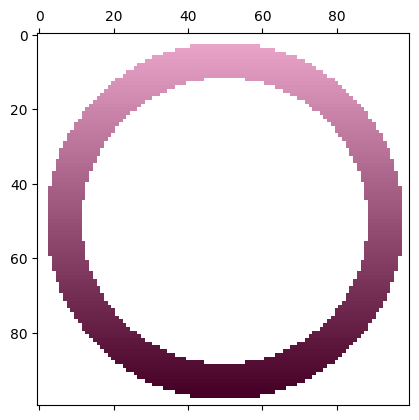} & 
\includegraphics[height = 0.15\linewidth]{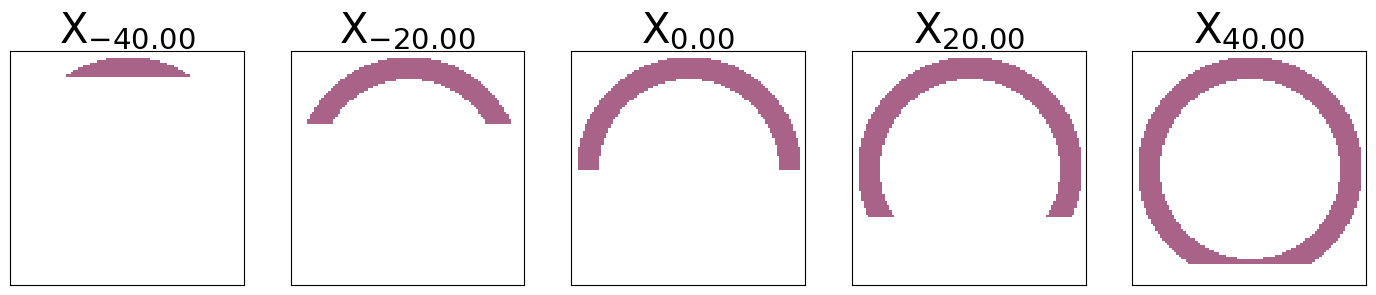} &
\includegraphics[height = 0.15\linewidth]{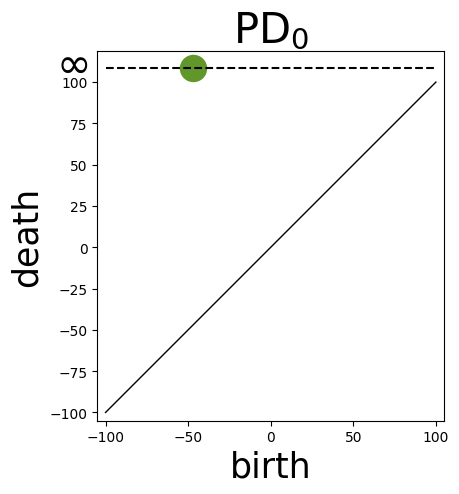} 
\includegraphics[height = 0.15\linewidth]{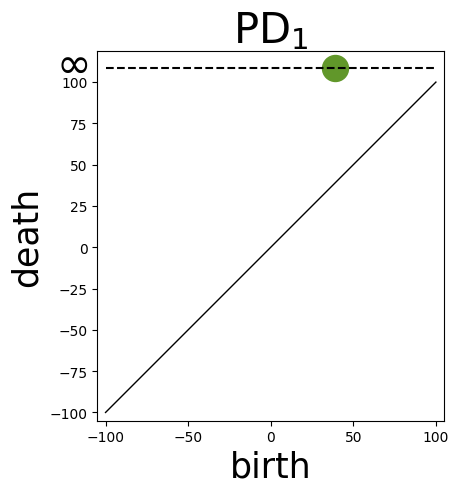} \\
\includegraphics[height = 0.15\linewidth]{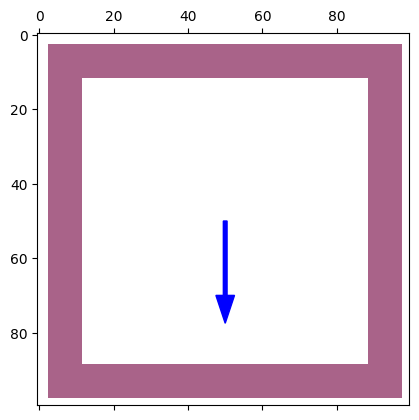} & 
\includegraphics[height = 0.15\linewidth]{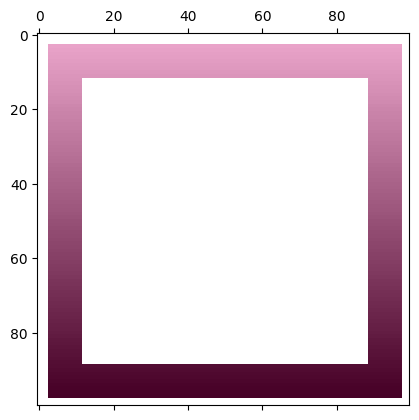} & 
\includegraphics[height = 0.15\linewidth]{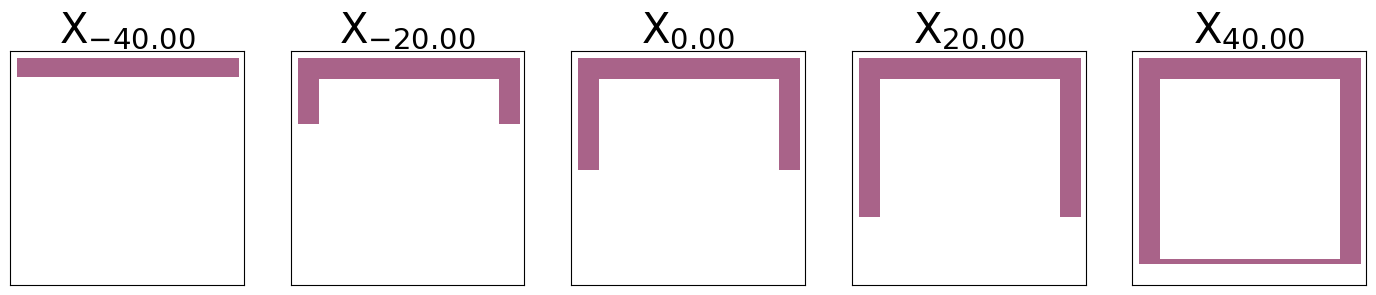} &
\includegraphics[height = 0.15\linewidth]{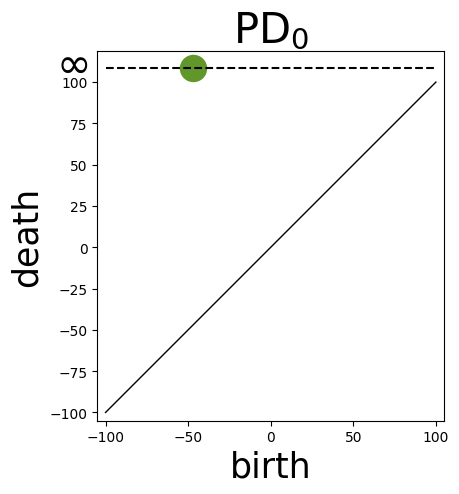} 
\includegraphics[height = 0.15\linewidth]{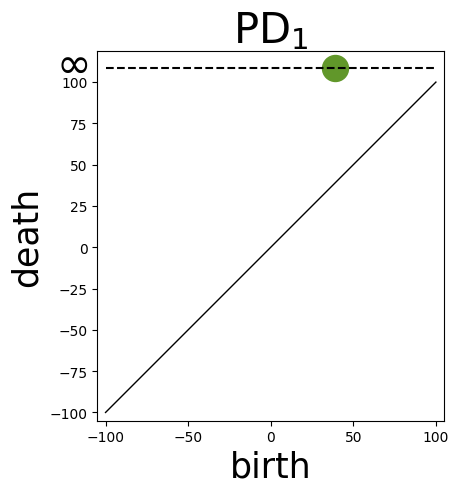} \\
\midrule

\includegraphics[height = 0.15\linewidth]{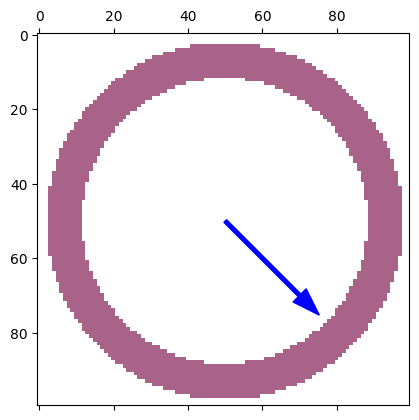} & 
\includegraphics[height = 0.15\linewidth]{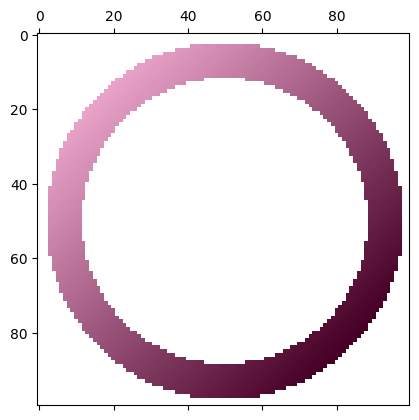} & 
\includegraphics[height = 0.15\linewidth]{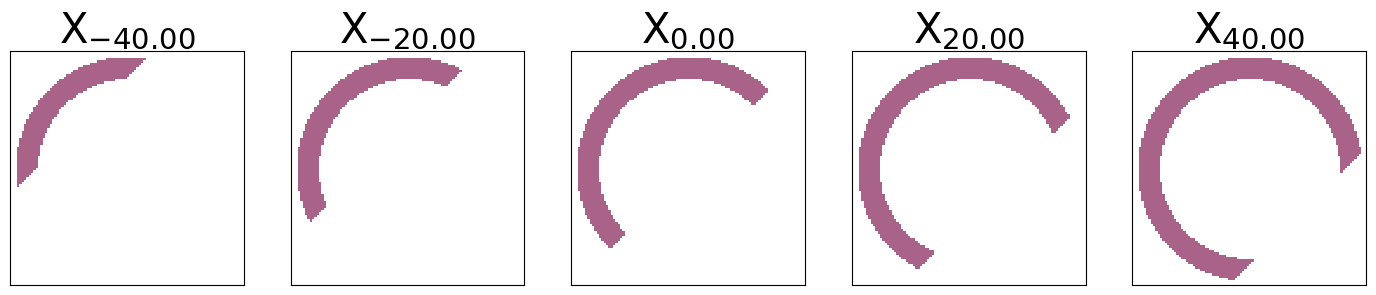} &
\includegraphics[height = 0.15\linewidth]{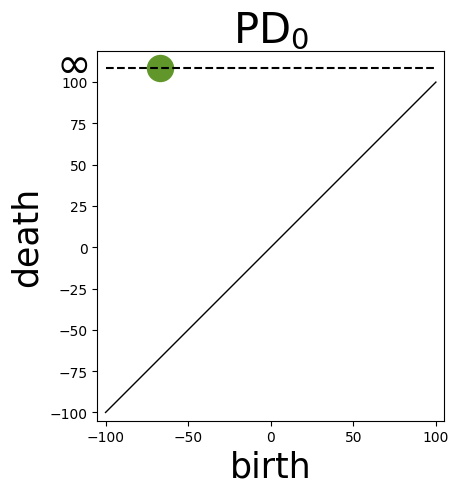} 
\includegraphics[height = 0.15\linewidth]{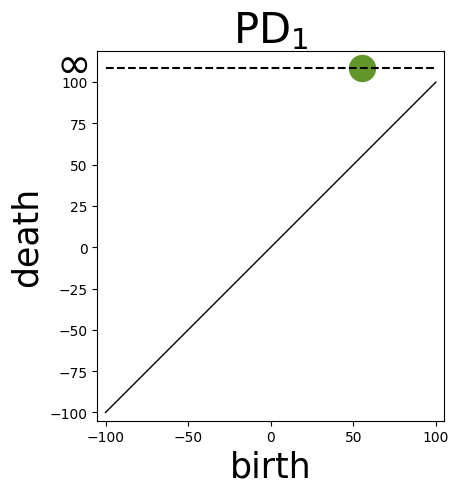} \\
\includegraphics[height = 0.15\linewidth]{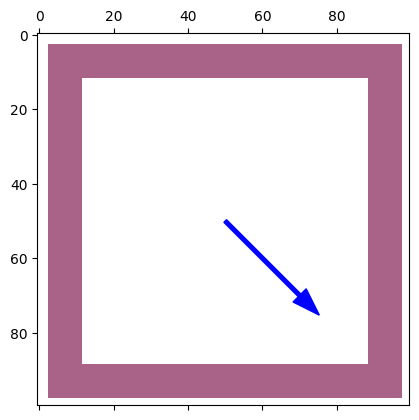} & 
\includegraphics[height = 0.15\linewidth]{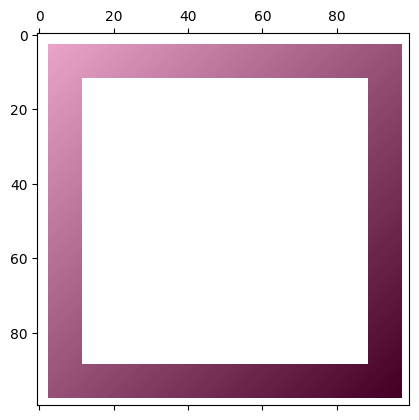} & 
\includegraphics[height = 0.15\linewidth]{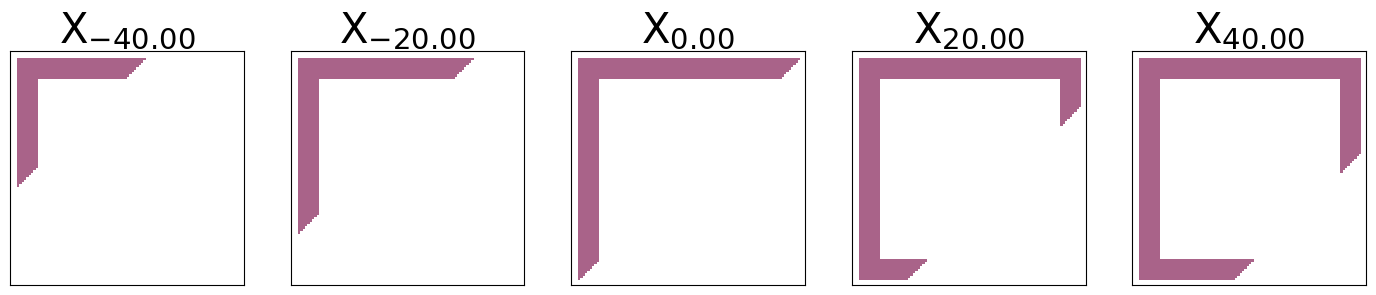} &
\includegraphics[height = 0.15\linewidth]{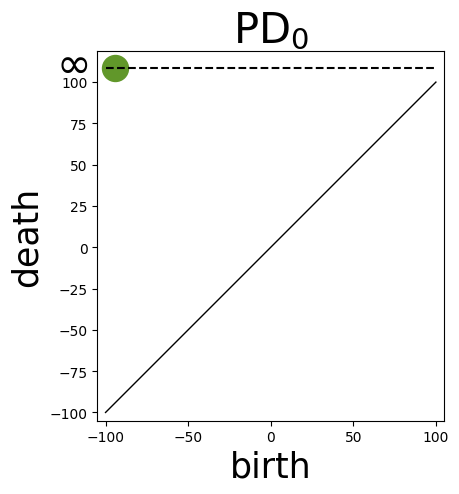} 
\includegraphics[height = 0.15\linewidth]{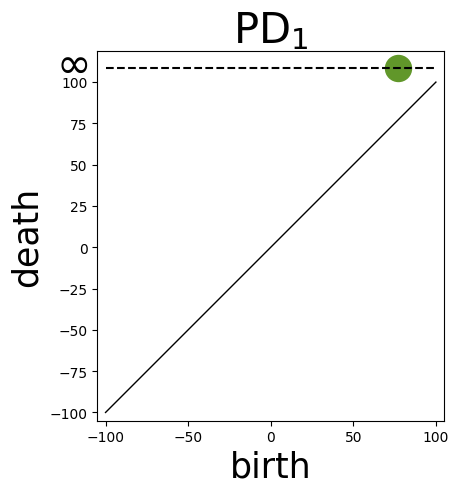} \\

\bottomrule
\end{tabular}
}
\caption{Illustration of the added value of $\GPHT$ over $\PH$ in every homological dimension: $\PH$ describes some topological and geometric features (such as the number, size or position of cycles) and therefore loses some information about the shape --- it is not injective, but considering multiple filtrations or ``directions" with $\GPHT$ can help to completely recover the shape. For the example circle and square, their $\PD$s are the same in any homological degree for some filtration (top rows), but considering different filtration can help to  discriminate between the two shapes (bottom rows).}
\label{fig:injectivity}
\end{figure}

Similarly to Example~\ref{E:networks}, in which one studies a tree by considering the radial filtration from a designated root note, one could more generally consider $\GPHT$ on the height, \dots, tubular or radial filtration with the parameter space $\mathbb{P}$ being a singleton: a particular, given hyperplane, \dots, line or a point (or more broadly, any single filtration, such as the greyscale filtration on images). 
While such an approach can be sufficient for particular applications, in general, it does not yield an injective descriptor of the shape. For instance, $\PH$ on the greyscale or height filtration from some direction might not be able to differentiate between the MNIST images of handwritten digits "0", "6" and "9" (see also Figure~\ref{fig:injectivity}). 
For probing shapes in practice, however, the ideal scenario lies somewhere in between ($\mathbb{P}$ being a singleton, which is not sufficient to capture the shape, and $\mathbb{P} = \AG(m, n)$ which is not possible in practice --- we cannot consider \emph{all} possible directions.) Ideally, we limit the $\GPHT$ to the set of all interesting flats (lines, planes, \dots, hyperplanes) that are sufficient to completely describe a shape.

\section{Continuity of $\DPHT(X)$} 
\label{sec:continuity}

In this section we work with affine Grassmannians and assume a choice of metric which is determined by principal angles.
We wish to show that $\GPHT(X)$ of an affine Grassmannian is continuous, so that a small perturbation to a flat will only cause a small perturbation in the diagrams.
The difficulty in doing this lies in converting between the affine Grassmannian metric and Euclidean notions of distance.
There is a notion of distance from a point to a plane in the doubly infinite affine Grassmannian by considering elements of $\R^n$ as $0$-dimensional flats.
However this metric is somewhat paradoxical from a Euclidean viewpoint, as all points are a strictly a positive distance away from a higher dimensional flat, even if they lie on the flat in the Euclidean sense.

\subsection{Some properties of singular values}
In this section we work with a metric on $\AG(m,n)$ determined by principal angles, such as Equation~\ref{eq:grassmann-distance}, where we recall that principal angles relate to singular values via Equation~\ref{eq:angle-singular}.

Recall that for a (real-valued) matrix $M$, the singular values of $M$ are the square roots of the eigenvalues of $M^TM$, where $M^TM$ is known to have a full decomposition of non-negative eigenvalues.
If $M$ is $m\times n$ then the singular value decomposition of $M$ is
\[
M = U\Sigma V^T
\]
where $U$ is an orthogonal $m\times m$ matrix, $V$ is an orthogonal $n\times n$ matrix, and $\Sigma$ is a diagonal $m\times n$ matrix whose non-zero entries are the singular values of $M$ in non-increasing order.
Efficient numerical methods for computing singular values (and principal angles) have been known since first being introduced by Jordan in 1875 \cite{jordan1875essai}.
One result which we will use in this section is from \cite[Theorem 1]{bjorck1973numerical}, which which shows that principal angles can be calculated algorithmically as well as through linear algebraic means.
Explicitly, suppose $\mathbb{A},\mathbb{B}\in\Gr(m,n)$ and let $M_\mathbb{A}$ and $M_\mathbb{B}$ be $m\times n$ matrices whose columns form an orthonormal basis of $\mathbb{A}$ and $\mathbb{B}$ respectively.
Then
\[
\cos\theta_1(\mathbb{A},\mathbb{B}) = \sigma_1(M_\mathbb{A}^TM_\mathbb{B}) = \max_{a\in\mathbb{A}}\max_{b\in\mathbb{B}} \{a\cdot b\,:\,\norm{a}=\norm{b}=1\}
\]
where $\sigma_1(M_\mathbb{A}^TM_\mathbb{B})$ is the largest singular value of $M_\mathbb{A}^TM_\mathbb{B}$.
Moreover, set $a_1$ and $b_1$ to be the unit vectors so that $a_1\cdot b_1=\sigma_1(M_\mathbb{A}^TM_\mathbb{B})$ and let $A_1=\mathrm{span}\{a_1\}$ and $B_1=\mathrm{span}\{b_1\}$.
Then we may determine the remaining principal angles by the following recursive procedure for $2\leq k\leq m$:
\begin{enumerate}[label=(\roman*)]
    \item We first note \[\cos\theta_{k}(\mathbb{A},\mathbb{B}) = \sigma_{k}(M_\mathbb{A}^TM_\mathbb{B}) = \max_{a\in\mathbb{A}\cap}\max_{b\in\mathbb{B}} \{a\cdot b\,:\,\norm{a}=\norm{b}=1, a\in A_{k-1}^\perp, b\in B_{k-1}^\perp\}\, .\]
    \item Set $a_k$ and $b_k$ to be the unit vectors so that $a_k\cdot b_k = \sigma_k(M_\mathbb{A}^TM_\mathbb{B})$.
    \item Set $A_k=\mathrm{span}\{a_1,\dots,a_k\}$ and $B_k=\mathrm{span}\{b_1,\dots,b_k\}$.
\end{enumerate}

\begin{theorem}\cite[Theorem 1]{bjorck1973numerical}
The above procedure returns the same principal angles as Equation~\ref{eq:angle-singular}.
\end{theorem}

In order to prove continuity of $\DPHT(X)$ we will require a classical result of Weyl on the perturbation theory of singular values \cite{weyl1912}.
In order to state this result, we will need to work with the following matrix norms.

\begin{definition}[Matrix norms]
Let $A$ be an $n\times k$ matrix.
The \emph{spectral norm} of $A$, denoted $\norm{A}_2$ is defined so that
\[
\norm{A}_2 = \max_{x\in\R^k:\norm{x}=1}\norm{Ax}.
\]
The \emph{Frobenius norm} of $A$, denoted $\norm{A}_F$ is the Euclidean norm when considering $A$ as an element of $\R^{nk}$. That is,
\[
\norm{A}_F^2 = \sum_{i=1}^n\sum_{j=1}^k(a_{ij})^2.
\]
\end{definition}

\begin{remark}
It is easy to show that $\norm{A}_2\leq\norm{A}_F$.
\end{remark}

\begin{lemma}[Weyl \cite{weyl1912}] \label{lem:weyl}
Suppose $A$ and $B$ are $n\times n$ matrices. Then for each $i=1,\dots,n$, 
\[
|\sigma_i(A)-\sigma_i(B)|\leq \norm{A-B}_2.
\]
\end{lemma}

\subsection{Convergence in an Affine Grassmannian implies Euclidean convergence}

In this section we show that convergence of affine flats of the same dimension in the affine Grassmannian implies convergence of the flats in a Euclidean sense.
The intuition is quite simple: a small perturbation of a flat should not disproportionately alter the directions it spans or its displacement from the origin.
We were unable to find any result for this in the literature so we include a complete proof here.
The problem is subtly difficult to prove, in particular when applying the embedding map $j:\AG(m,n)\to\Gr(m+1,n+1)$ while using a large number of variables.

\begin{lemma} \label{lem:sing-val-converge}
Suppose $(P_k)$ is a sequence of flats in $\AG(m,n)$ such that $P_k\to P$.
Let $M_{P_k}$ and $M_P$ be matrices whose columns are orthonormal bases of $j(P_k)$ and $j(P)$ respectively. Then the singular values of $(M_{P_k})^TM_P$ all converge to $1$ independent of the choice of matrices.
\end{lemma}
\begin{proof}
$P_k\to P$ in the metric if and only if all principle angles converge to zero.
Principal angles are related to singular values via the equation $\cos\left(\theta_i(P_k,P)\right)=\sigma_i\left((M_{P_k})^TM_P\right)$, which must converge to $1$.
The convergence is independent of the matrix representatives because singular values are fixed under coordinate transformation.
\end{proof}

The key ingredient for our proof will be to apply Lemma~\ref{lem:weyl} to uniformly bound a singular value of $(M_{P_k})^TM_P$ below 1.
In order to achieve this, we will need the following bound.

\begin{lemma} \label{lem:bound}
Suppose $x$ and $p$ are two vectors such that $\norm{x-p}\geq\Delta>0$.
Then there exists a positive bound $B_{p,\Delta}<1$ constant with respect to $x$ such that 
\[
\left|\frac{1+p\cdot x}{\sqrt{(1+\norm{p}^2)(1+\norm{x}^2)}}\right|\leq B_{p,\Delta}
\]
\end{lemma}
We prove Lemma~\ref{lem:bound} in Appendix~\ref{app:bound-proof} since it is rather technical and does not provide notable insight of Affine Grassmannians.

Finally, we will also need to use the fact that the deaffine map is continuous.
When considering $\AG(m,n)$ as a vector bundle over $\Gr(m,n)$, the deaffine map is continuous by definition.
However, the metric on $\AG(m,n)$ is not inherited from this model, so we include a brief proof here for completeness.

\begin{lemma}
The deaffine map $\pi:\AG(m,n)\to\Gr(m,n)$ is continuous.
\end{lemma}
\begin{proof}

Let $B_{\epsilon}(\mathbb{A})$ be an open ball in $\Gr(m,n)$. Its preimage under $\pi$ is given by 
\[
\left \{\mathbb{B}+b \mid \mathbb{B}\in B_{\epsilon}(\mathbb{A}),\;\; b\in \mathbb{R}^n\right \} =\bigcup_{b\in \mathbb{R}^n} \left \{\mathbb{B}+b \mid \mathbb{B}\in B_{\epsilon}(\mathbb{A})\right \}\, .
\]
Each set $\left \{\mathbb{B}+b \mid \mathbb{B}\in B_{\epsilon}(\mathbb{A})\right \}$ is open in the subspace topology induced on $\AG(m,n)$ by the embedding $\AG(m,n)\hookrightarrow \Gr(m+1,n+1)$, since   $\left \{\mathbb{B}+b \mid \mathbb{B}\in B_{\epsilon}(\mathbb{A})\right \}$ may be identified with an open ball in $ \Gr(m+1,n+1)$ centered at $\mathrm{span}\{B',b'+e_{n+1}\}$, where $B$ is a choice of basis vectors for $\mathbb{B}$, and we recall that the notation $B'$ means that we are identifying $B\in \mathbb{R}^{n\times m}$ with an element of $\mathbb{R}^{(n+1)\times (m+1)}$ by adding zeros in the additional components. 
Thus, the union of all such open sets is open in $\AG(m,n)$.
\end{proof}

\begin{proposition} \label{prop:ag-implies-euc}
Suppose $(P_k)$ is a sequence of flats in $\AG(m,n)$ such that $P_k\to P$.
Then $P_k\to P$ in the Euclidean sense.
In particular, if $P_k = \mathbb{A}_k+b_k$ and $P=\mathbb{A}+b$ then
\begin{enumerate}
    \item For every $k$ there exist orthonormal bases $\{u_{1,k},\dots,u_{m,k}\}$ and $\{v_{1,k},\dots,v_{m,k}\}$ of $\mathbb{A}_k$ and $\mathbb{A}$ (respectively) such that $\norm{u_{i,k}-v_{i,k}}\to 0$ for each $i$.
    \item $\norm{b_k-b}\to 0$.
\end{enumerate}
\end{proposition}
\begin{proof}
If $m=0$ then the first statement is vacuously true.
Otherwise, the deaffine map $\pi:\AG(m,n)\to\Gr(m,n)$ is continuous, so $P_k\to P$ in the Affine Grassmannian implies $\mathbb{A}_k\to \mathbb{A}$ in the Grassmannian.
In particular, the principal angles between $\mathbb{A}_k$ and $\mathbb{A}$ converge to zero.
Calculating principle angles via the algorithmic method yields $\{u_{1,k},\dots,u_{m,k}\}\subseteq \mathbb{A}_k$ and $\{v_{1,k},\dots,v_{m,k}\}\subseteq \mathbb{A}$ such that for each $i=1,\dots,m$ we have $\cos(\theta_i)=u_{i,k}\cdot v_{i,k}$, $\norm{u_{i,k}}=\norm{v_{i,k}}=1$, $u_{i,k}\cdot u_{j,k}=v_{i,k}\cdot v_{j,k} = 0$ for each $j<i$.
In particular $\{u_{1,k},\dots,u_{m,k}\}$ and $\{v_{1,k},\dots,v_{m,k}\}$ are orthonormal bases by construction and $\theta_i\to 0$ implies $u_{i,k}\cdot v_{i,k}\to 1$ for each $i$.
Since
\[
\norm{u_{i,k}-v_{i,k}}^2 = (u_{i,k}-v_{i,k})\cdot(u_{i,k}-v_{i,k}) = \norm{u_{i,k}}^2+\norm{v_{i,k}}^2-2u_{i,k}\cdot v_{i,k} = 2-2u_{i,k}\cdot v_{i,k}\to 0
\]
this proves the first statement.

For the second statement, suppose towards a contradiction that $(b_k)$ does not converge to $b$ in the Euclidean sense.
Then passing to a subsequence if necessary, we may assume there exists $\Delta>0$ such that $\norm{b_k-b}\geq\Delta$ for each $k\in\N$.
To apply this assumption, we must obtain an expression for the convergence of $P_k\to P$.
For this reason, let $\{u_{1,k},\dots,u_{m,k}\}$ and $\{v_{1,k},\dots,v_{m,k}\}$ be as above.
The embedding map $j:\AG(m,n)\to\Gr(m+1,n+1)$ is defined explicitly on $P_k$ and $P$ so that
\[
j(P_k) = \mathrm{span}\left\{u'_{1,k},\dots,u'_{m,k},\frac{b_k+e_{n+1}}{\sqrt{1+\norm{b_k}^2}}\right\}
\] and 
\[j(P) = \mathrm{span}\left\{v'_{1,k},\dots,v'_{m,k},\frac{b+e_{n+1}}{\sqrt{1+\norm{b}^2}}\right\}.
\]

In particular, the two expressions span over orthonormal bases of $j(P_k)$ and $j(P)$, so let $M_k$ and $N_k$ be the matrices whose columns are these basis vectors.
By Lemma~\ref{lem:sing-val-converge}, $P_k\to P$ implies $\sigma_i((M_k)^TN_k)\to 1$ for each $i$, so we construct $(M_k)^TN_k$ explicitly.
Let $\Tilde{M}_k$ and $\Tilde{N}_k$ be the top left $m\times m$ minors of $M_k$ and $N_k$, whose columns are explicitly the orthonormal bases of $\mathbb{A}_k$ and $\mathbb{A}$.
Then $(M_k)^TN_k$ is of the form
\[
\begin{pmatrix}
 & (\tilde{M}_k)^T\tilde{N}_k & & \vline & \frac{(\tilde{M}_k)^T(b+e_{n+1})}{\sqrt{1+\norm{b}^2}}\\
    \hline
 & \frac{(b_k+e_{n+1})^T\tilde{N}_k}{\sqrt{1+\norm{b_k}^2}} & & \vline & \frac{1+b\cdot b_k}{\sqrt{(1+\norm{b}^2)(1+\norm{b_k}^2)}}\\
\end{pmatrix}
\]
The top left component is an $m\times m$ matrix whose $(i,j)^\mathrm{th}$ term is exactly $u_{i,k}\cdot v_{j,k}$.
When $i=j$ this converges to $1$ by the first statement.
When $i\neq j$ we have 
\[
|u_{i,k}\cdot v_{j,k}| = |u_{i,k}\cdot (v_{j,k}-u_{j,k)})|\leq \norm{v_{j,k}-u_{j,k)}}\to 0.
\]
Thus $(\tilde{M}_k)^T\tilde{N}_k$ converges pointwise to the identity matrix $I_m$, and therefore converges in the Frobenius norm to $I_m$.

The bottom left component of $(M_k)^TN_k$ is an $m$-dimensional row vector whose $i^\mathrm{th}$ entry is exactly
\[
\frac{(b_k+e_{n+1})\cdot v_{i,k}}{\sqrt{1+\norm{b_k}^2}} = \frac{b_k\cdot v_{i,k}}{\sqrt{1+\norm{b_k}^2}}.
\]
And since $b_k$ is orthogonal to $\mathbb{A}_k$ we have
\[
|b_k\cdot v_{i,k}| = |b_k\cdot (v_{i,k}-u_{i,k})|\leq \norm{b_k}\norm{v_{i,k}-u_{i,k}}
\]
So 
\[
\left|\frac{(p_k+e_{n+1})\cdot v_{i,k}}{\sqrt{1+\norm{p_k}^2}}\right| \leq \frac{\norm{p_k}\norm{v_{i,k}-u_{i,k}}}{\sqrt{1+\norm{p_k}^2}} \leq \norm{v_{i,k}-u_{i,k}}\to 0.
\]
Therefore the bottom left component of $(M_k)^TN_k$ converges to zero pointwise, and therefore in norm.
A symmetric argument shows that the top right component is a column vector also converging to the zero vector in norm.

Finally, the bottom right component of $(M_k)^TN_k$ is a constant whose absolute value by Lemma~\ref{lem:bound} is bounded by $B_{b,\Delta}<1$ independent of $b_k$.

To complete the proof, we set $\tilde{I}$ to be the matrix of the form
\[
\begin{pmatrix}
 & I_m & \vline & 0\\
    \hline
 & 0 & \vline & 0\\
\end{pmatrix}
\]
and use Lemma~\ref{lem:weyl} to bound $\sigma_{m+1}((M_k)^TN_k)$.
Note that $\sigma_{m+1}(\tilde{I})=0$, so 
\begin{align*}
\sigma_{m+1}((M_k)^TN_k)^2& \leq \norm{(M_k)^TN_k-\tilde{I}}^2_2\\
& \leq \norm{(M_k)^TN_k-\tilde{I}}^2_F \\
& = \norm{(\tilde{M}_k)^T\tilde{N}_k-I_m}^2_F + \norm{\frac{(b_k)^T\tilde{N}_k}{1+\norm{b_k}^2}}^2 + \norm{\frac{(\Tilde{M}_k)^Tb}{1+\norm{b}^2}}^2 \\
& \;\;\;\;\;\; + \left|\frac{1+b\cdot b_k}{\sqrt{(1+\norm{b}^2)(1+\norm{b_k}^2)}}\right|^2
\end{align*}
The first three terms converge to 0 and the final term is bounded above by $B_{b,\Delta}^2$, so for $k$ sufficiently large we must have 
\[
|\sigma_{m+1}((M_k)^TN_k)| \leq \sqrt{\frac{1 + B_{b,\Delta}^2}{2}}<1.
\]
But $\sigma_{m+1}((M_k)^TN_k)\to 1$, so we arrive at a contradiction.
\end{proof}

\subsection{Bottleneck continuity theorem}

Now that we can guarantee convergence in an affine Grassmannian metric implies convergence in a Euclidean sense, our proof of continuity of the affine Grassmannian $\DPHT$ with respect to bottleneck distances of persistence diagrams is very similar to the proof of the classical $\PHT$ in \cite{TMB2014}.
In particular, we will apply bottleneck stability of persistence diagrams and then show that the supremum distance between two ``distance to flat'' functions is uniformly bounded.
Where our method diverges from \cite{TMB2014} is that we cannot use a Lipschitz continuity argument, and instead apply a sequential continuity argument using Euclidean convergence results from Proposition~\ref{prop:ag-implies-euc}.

\begin{theorem}[$\GPHT$ is continuous with respect to bottleneck distance] \label{thm:bottleneck-continuity}
Let $X$ be a constructible set, $n\in\N$ and $0\leq m<n$.
Then $\DPHT(X):\AG(m,n)\to \mathcal{D}^n$ is continuous with respect to bottleneck distance.
\end{theorem}
The assumption that $X$ is constructible is imposed to be consistent with what we require in Section~\ref{sec:injectivity}, while also applying bottleneck stability (Theorem~\ref{thm:bottleneck-stability}).
In the following proof, $X$ being constructible can be weakened to definable and bounded.
We require $X$ to be bounded in order to control bounds on the supremum metric of continuous functions on $X$, and we require that $X$ be definable and that distance-to-flat-functions are continuous in order to apply bottleneck stability.

\begin{lemma}\label{L:d_P lipsch}
Suppose $X\subset\R^n$ is a constructible set, $P\in\AG(m,n)$ and $f_P:X\to[0,\infty)$ is the distance-to-flat function.
Then $f_P$ is  $(m+1)$-Lipschitz continuous.
\end{lemma}
\begin{proof}

We recall from elementary linear algebra that distances from a point to a plane are measured by orthogonal projection.

Let $P=\mathbb{A}+b$ with $\{v_1,\dots,v_m\}$ an orthonormal basis of $\mathbb{A}$ and $b$ be the displacement of $P$ from the origin.
Then the orthogonal projection of $x$ onto $P$ is
\[
x^\perp = b + \sum_{i=1}^m((x-p)\cdot v_i)v_i,
\]
where the sum is empty when $m=0$.
Further, the vector $b$ is orthogonal to the underlying vector space by construction, so in particular
\begin{align*}
f_P(x) = \norm{x-b - \sum_{i=1}^m(x\cdot v_i)v_i}.
\end{align*}
Then
\begin{align*}
|f_P(x)-f_P(y)| & = \left|\norm{x-b - \sum_{i=1}^m(x\cdot v_i)v_i}-\norm{y-b - \sum_{i=1}^m(y\cdot v_i)v_i}\right| \\
& \leq \norm{x-y - \sum_{i=1}^m((x-y)\cdot v_i)v_i} \\
& \leq \norm{x-y} + \sum_{i=1}^m\norm{x-y}\norm{v_i}^2 = (m+1)\norm{x-y}
\end{align*}
which gives $(m+1)$-Lipschitz continuity.
\end{proof}

\begin{proof}[Proof of Theorem~\ref{thm:bottleneck-continuity}]
Let $P_k=\mathbb{A}_k+b_k$ and $P=\mathbb{A}+b$.
As $f_P$ and $f_{P_k}$ are continuous and $X$ is constructible, we may apply bottleneck stability \cite{CSEH06} of persistence diagrams.
This implies we need only show the map $P\mapsto f_P$ is continuous with respect to the supremum metric.

Let $k\in \N$ and let $\{u_{1,k},\dots,u_{m,k}\}$ and $\{v_{1,k},\dots,v_{m,k}\}$ be orthonormal bases $\mathbb{A}_k$ and $\mathbb{A}$ respectively as described in Proposition~\ref{prop:ag-implies-euc}.
We need to show that $f_{P_k}\to f_P$ in the supremum metric.
With distances described as in the previous Lemma, this gives
\begin{align}
|f_P(x)-f_{P_k}(x)| & = \left|\norm{x-b - \sum_{i=1}^m(x\cdot v_{i,k})v_{i,k}}-\norm{x-b_k - \sum_{i=1}^m(x\cdot u_{i,k})u_{i,k}}\right|\\
& \leq \norm{b_k-b + \sum_{i=1}^m\left((x\cdot u_{i,k})u_{i,k}-(x\cdot v_{i,k})v_{i,k}\right)}
\end{align}
Note that 
\[
(x\cdot u_{i,k})u_{i,k}-(x\cdot v_{i,k})v_{i,k} = (x\cdot u_{i,k})(u_{i,k}-v_{i,k}) - (x\cdot (v_{i,k}-u_{i,k}))v_{i,k}
\]
which means
\begin{align*}
|f_P(x)-f_{P_k}(x)| & \leq \norm{b_k-b}+ \sum_{i=1}^m\left(\norm{x}\norm{u_{i,k}}\norm{u_{i,k}-v_{i,k}} + \norm{x}\norm{v_{i,k}-u_{i,k}}\norm{v_{i,k}}\right) \\
& \leq \norm{b_k-b}+ 2\norm{x}\sum_{i=1}^m\norm{u_{i,k}-v_{i,k}}
\end{align*}
Finally, since $X$ is compact, $\norm{x}$ is uniformly bounded and the remaining terms converge to zero by Proposition~\ref{prop:ag-implies-euc} independent of $x$.
Therefore $f_{P_k}\to f_P$ in the supremum metric and we are done.
\end{proof}

\begin{remark}
While the distance-to-flat functions are Lipschitz and the distance-to-flat $\PHT$s are sublevel set filtrations of these functions, $\DPHT(X):\AG(m,n)\to\mathcal{D}^n$ will not be Lipschitz except for trivial examples.
This in contrast to the classical $\PHT$, which is 1-Lipschitz continuous with a fairly elementary proof \cite[Lemma~2.1]{TMB2014}.
The fundamental barrier is that most metrics of $\AG(m,n)$ are bounded when defined in terms of principle angles.
For instance, if $X$ is the origin in $\R^2$ and $\ell_r$ is the horizontal line $y=r$, then $\theta(\ell_0,\ell_r)=\arccos((1+r^2)^{-1/2})$ which approaches $\pi/2$ as $r\to\infty$ but $\norm{f_{\ell_0}-f_{\ell_r}}_\infty = r$ is unbounded. 
\end{remark}

\subsection{Wasserstein continuity theorem}
A main ingredient in the proof of continuity of $\PHT(X)$ with respect to the Wasserstein distance on $\D^n$
is to show that the total persistence (i.e., the sum of all the lifespans over all points in a diagram) of each persistence diagram is bounded. 
 This is done in Lemma 3.4 of  \cite{CST22}, which guarantees that the homeomorphism type of each sublevel set   changes only finitely many times across the filtration parameter. The proof of Lemma 3.4 carries over to our setting in a straightforward way. For completeness, we include here the corresponding statement in our setting:
 
 \begin{lemma}\label{L:finitely many changes}
For any constructible set $X\subset \mathbb{R}^n$ and any flat $P\in \AG(m,n)$ the homeomorphism type of the  sublevel set $X^{-}_r(\dist_P)=\{x\in X \mid \dist(x,P)\leq r\}$   changes only finitely many times as a function of $r\in \mathbb{R}$. Furthermore, there is a bound $B_X\geq0$, which depends on $X$, on the total number of changes in homeomorphism type over all flats $P$. 
\end{lemma}

\begin{proof}
The proof of Lemma 3.4 from \cite{CST22} carries over almost verbatim to our setting, with the main difference being that we consider $\mathbb{P}=\AG(m,n)$ instead of $\mathbb{P}=S^{n-1}$, and $d_P$ instead of $h_v$. We have showed that $\AG(m,n)$ is definable in Example \ref{E: aff grass def}, and that $d_P$ is definable in Example \ref{E: distance def}.
\end{proof}

\begin{theorem}[$\DPHT$ is continuous with respect to Wasserstein distance]
Let $X\subset \mathbb{R}^n$ be a  constructible subset, $n\in\N$ and $0\leq m<n$.
Then $\PHT(X):\AG(m,n)\to \mathcal{D}^n$ is continuous with respect to $p$-Wasserstein distance for $p> 1$.
\end{theorem}

\begin{proof}
By Lemma \ref{L:finitely many changes} we have that the total persistence of each persistence diagram is bounded. Furthermore, we have shown that the distance-from-flat functions are Lipschitz in Lemma \ref{L:d_P lipsch}, and that the function $P\mapsto d_P$ is continuous with respect to the supremum metric in the Proof of Theorem \ref{thm:bottleneck-continuity}. Finally, Wasserstein continuity for $p> 1$ follows from  \cite[Wasserstein Stability Theorem]{wasserstein-stability}. 

\end{proof}

\subsection{Instability of $\DPHT$}\label{remark:instability}

\begin{figure}[h]
\centering
\begin{tabular}{c|c|c|c}
\toprule
Shape $X$ & $d_P$ & Distance-from-flat filtration & $\PH(X)$ \\
\midrule

\includegraphics[height = 0.1\linewidth]{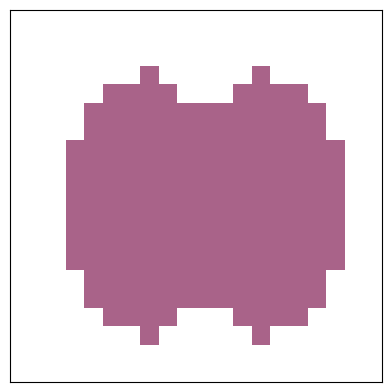} & 
\includegraphics[height = 0.1\linewidth]{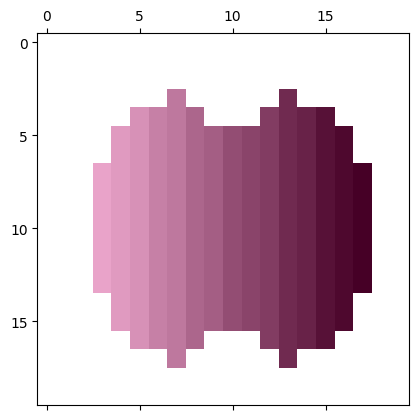} &
\includegraphics[height = 0.1\linewidth]{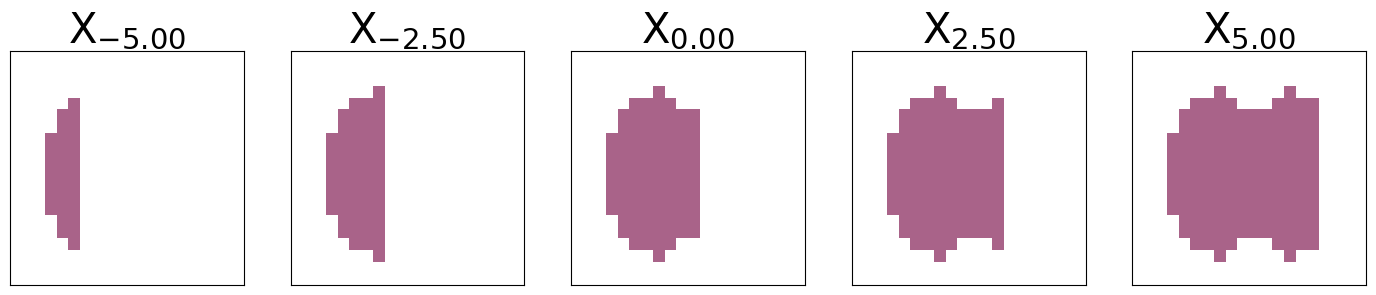} &  
\includegraphics[height = 0.1\linewidth]{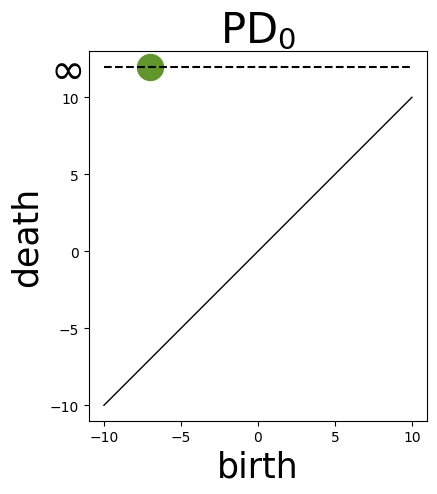} 
\includegraphics[height = 0.1\linewidth]{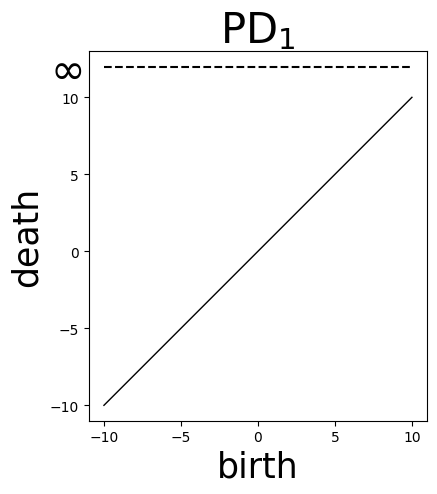} \\

\includegraphics[height = 0.1\linewidth]{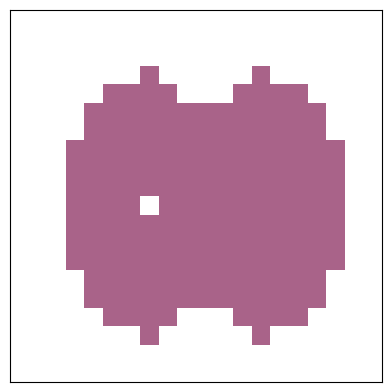} & 
\includegraphics[height = 0.1\linewidth]{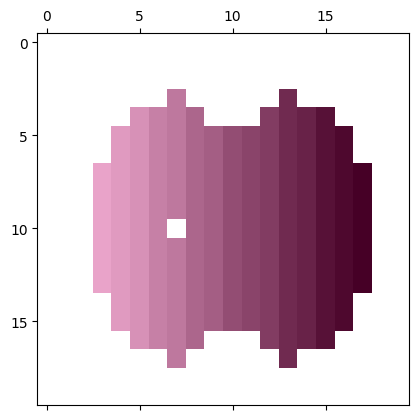} &
\includegraphics[height = 0.1\linewidth]{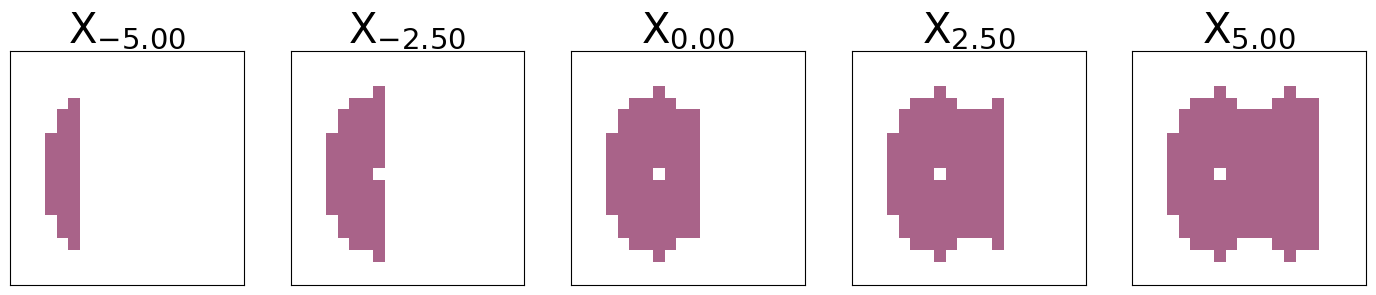} & 
\includegraphics[height = 0.1\linewidth]{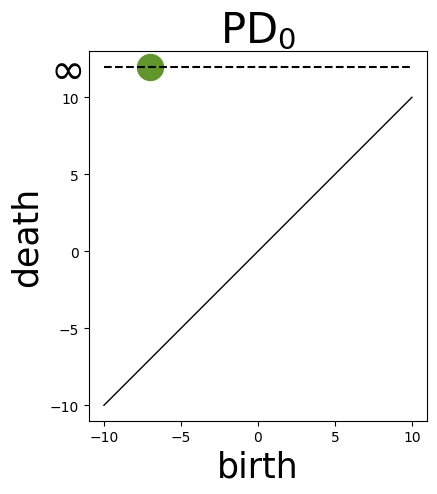} 
\includegraphics[height = 0.1\linewidth]{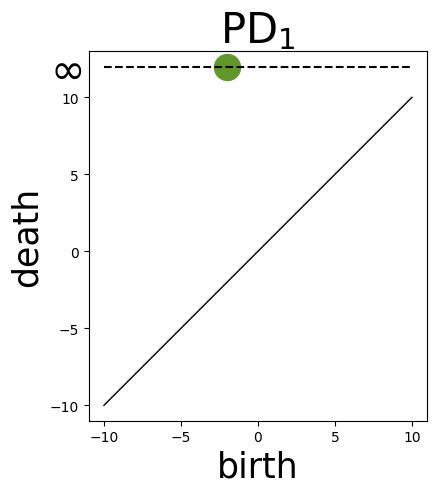} \\

\bottomrule
\end{tabular}
\caption{Distance-from-flat $\DPHT$ is not stable, since a small shape perturbation can yield a large change in $\PHT$. For the example shapes $X$ and $Y$ which are at small distance $d(X, Y),$ and any flat $P \in \AG(m, 2),$ we have that $d(\PD_1(X, d_P), \PD_1(Y, d_P)) = \infty$, since even a single-pixel hole results in an essential $1$-dimensional cycle.}
\label{fig:stability}
\end{figure}

We emphasise that while $\DPHT(X)$ is continuous --- namely, if $P, P' \in \AG(m, n)$ are close, then $\PD_k(X, d_P)$ and $\PD_k(X, d_{P'})$ are close ---, the function 
\[
\DPHT\colon CS(\mathbb{R}^n) \to  \mathrm{Fun}\left(\AG(m,n),\mathcal{D}^n\right ) 
\] is \emph{not} continuous: if shapes $X$ and $X'$ are close, it is not always true that $\PD_k(X, d_P)$ and $\PD_k(X', d_{P})$ are close. Commonly, in the TDA literature, we would say that $\DPHT$ is not stable. 
The issue arises whenever $X$ and $X'$ have different homology, or Betti numbers, i.e., their sublevel-set $\PD$s have a different number of essential classes (independently of how close $X$ and $X'$ are in, e.g., the Gromov-Hausdorff distance), since in this case $d(\PD_k(X, d_P), \PD_k(X', d_{P})) = \infty,$ for $d$ being the bottleneck or Wasserstein metric (Figure~\ref{fig:stability}). We note that this is a particular feature of sublevel-set filtrations; persistent homology of  geometric complexes such as Vietoris-Rips or \v{C}ech complexes is known to be stable on the set of totally bounded  metric spaces together with the Gromov-Hausdorff distance, see \cite[Theorem 5.2]{CdSO14}.

Luckily, recent work \cite{turner2024extended} on the extended persistent homology transform for the classical $\PHT$ rectifies this problem, albeit by restricting the class of allowed spaces to triangulated manifolds with boundary. We leave the question of exploring possible stabilisations of our  generalised $\PHT$s for future work. We end this Section with Example \ref{E:volcano island}, which  illustrates a possible application scenario in which small perturbations with respect to which  $\PHT$s are unstable  may arise.

\begin{figure}[h]
    \centering
    \includegraphics[width=0.49\linewidth]{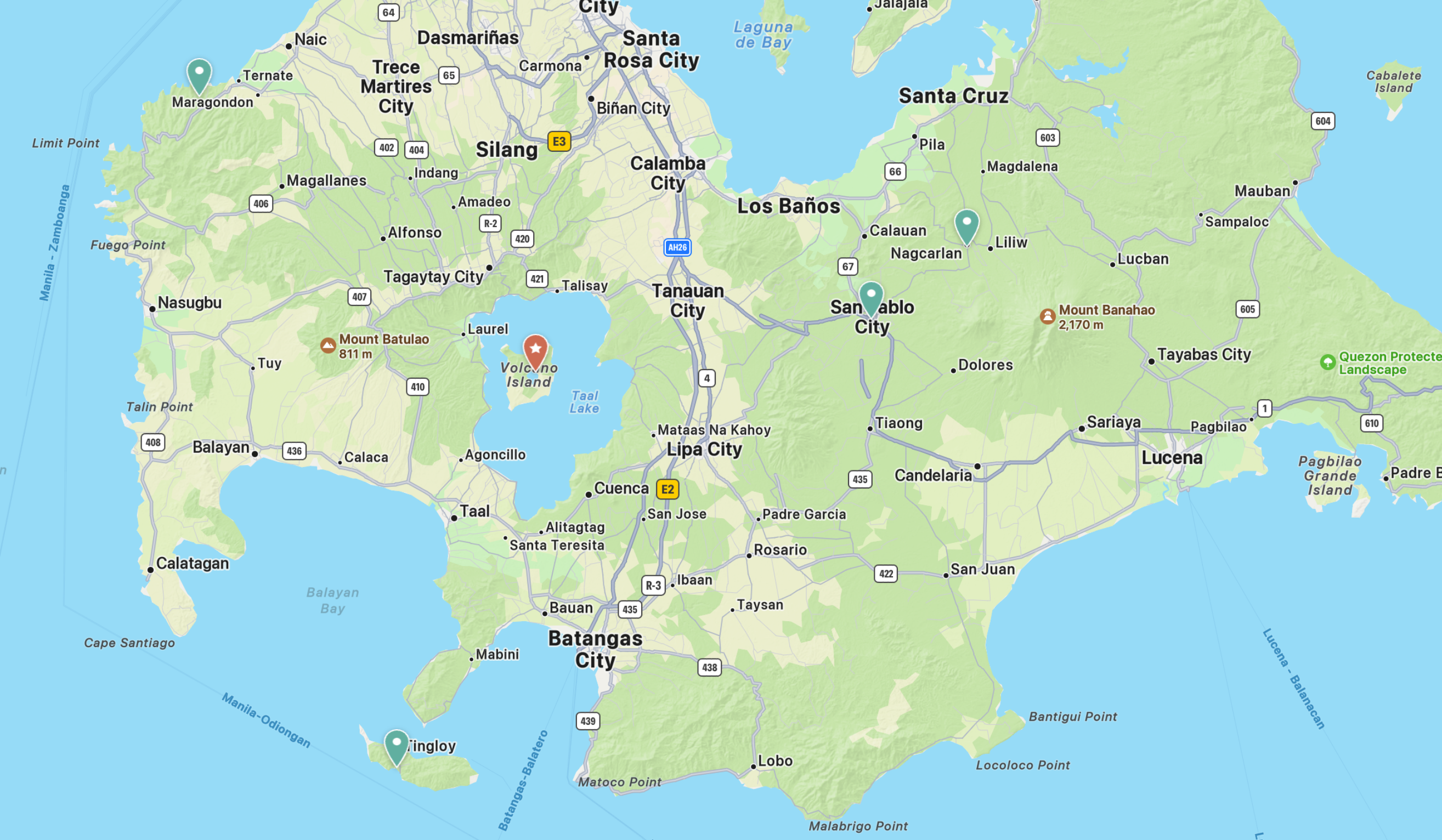}
    \includegraphics[width=0.49\linewidth]{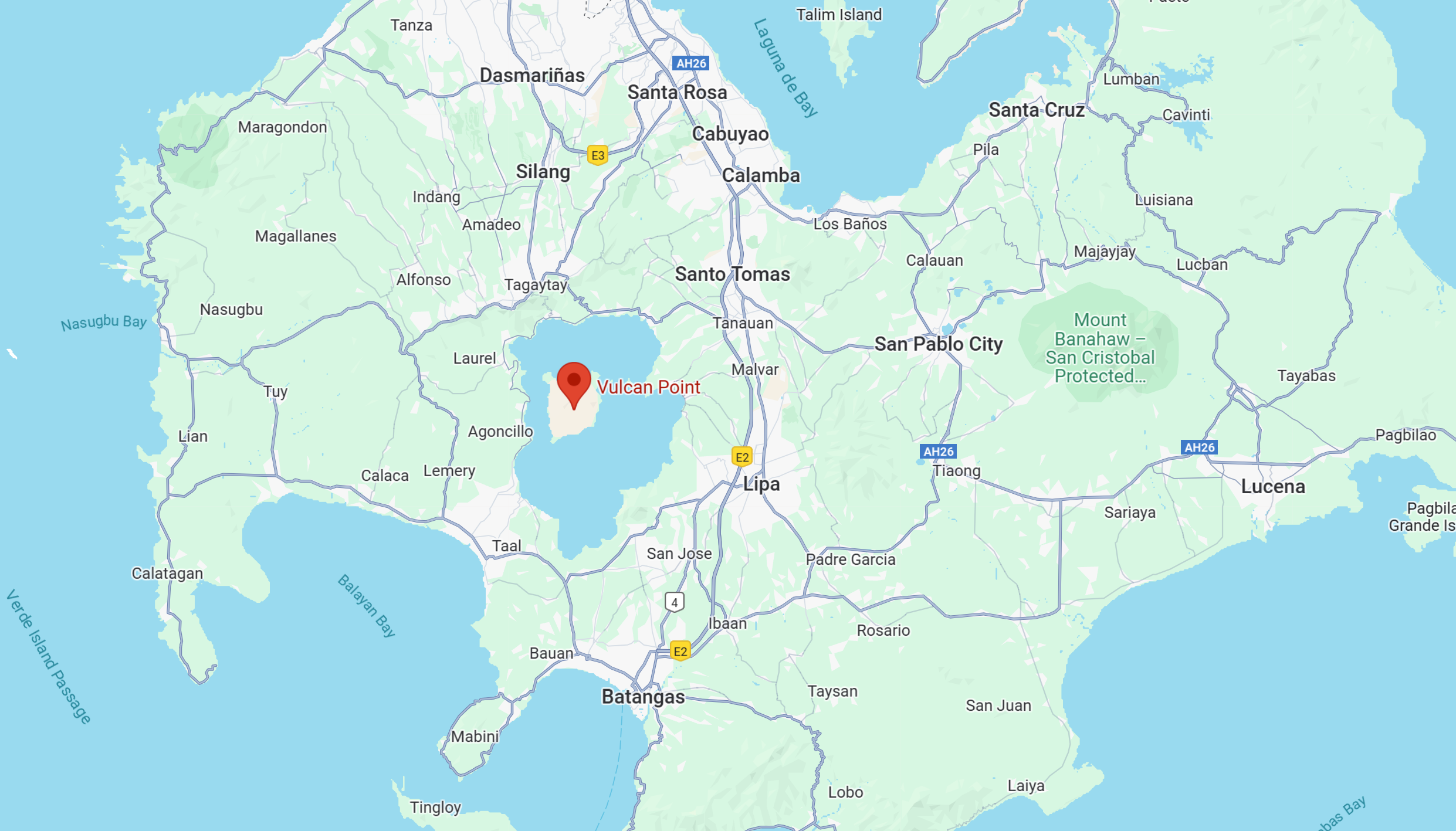}
    \includegraphics[width=0.49\linewidth]{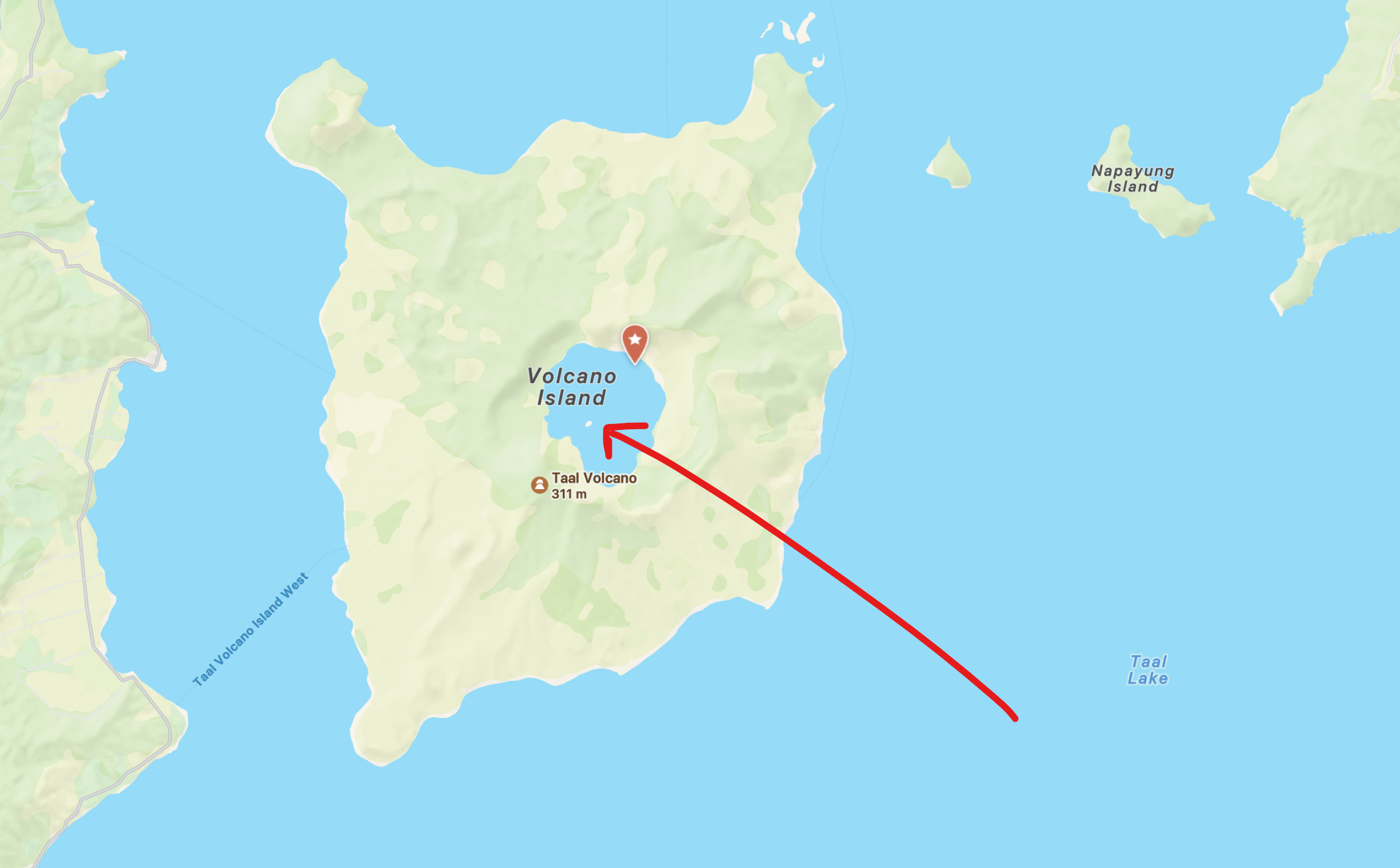}
    \includegraphics[width=0.49\linewidth]{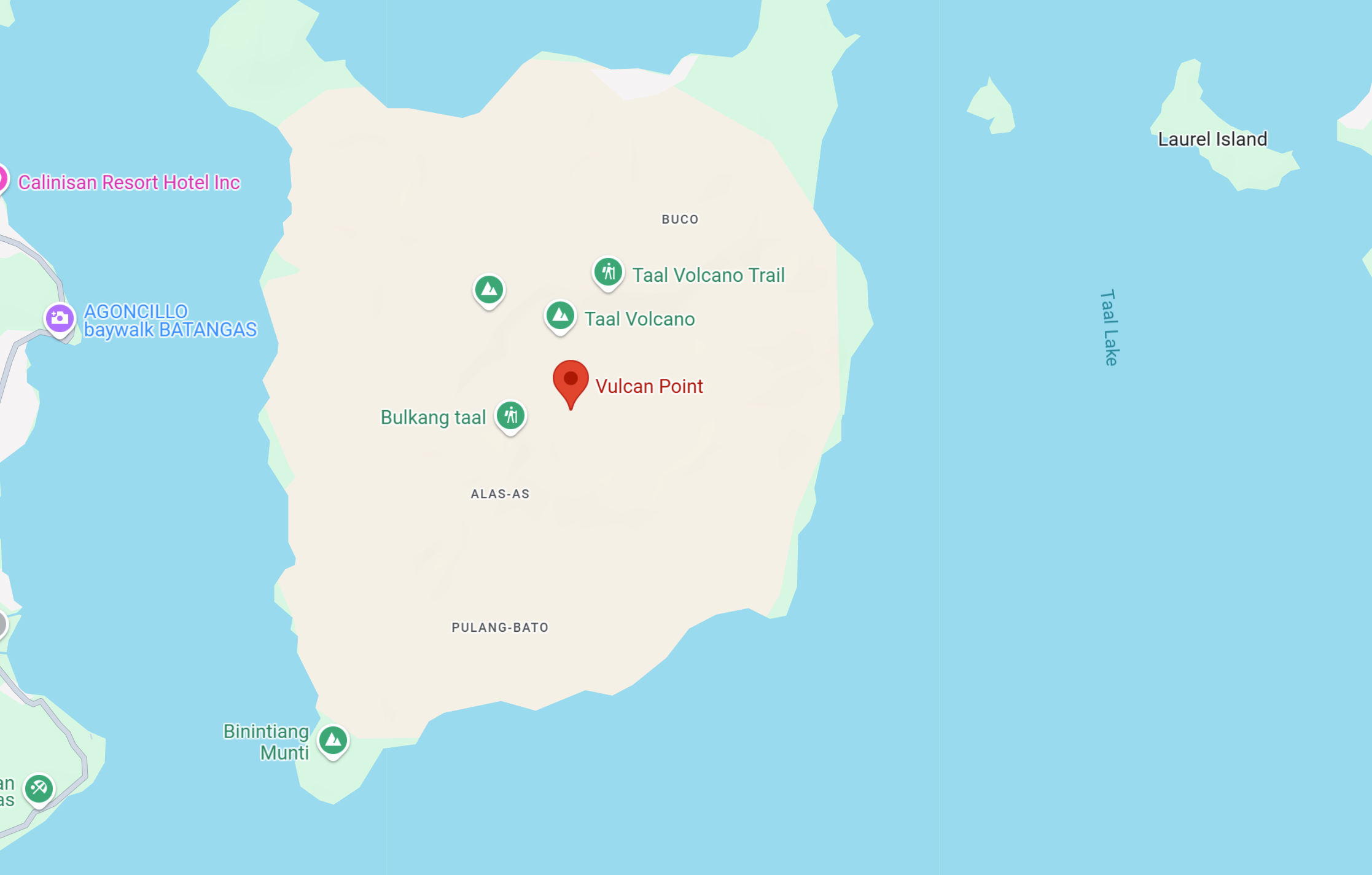}
    \caption{Top: maps of part of the island of Luzon containing Vulcan Point as viewed on Apple Maps (left \cite{apple}) and Google Maps (right \cite{google}). Bottom: maps of part of Taal Lake as viewed on Apple Maps (left, with an arrow pointing at Vulcan Point) and Google Maps (right, which is better up-to-date and no longer indicates the presence of  Taal Crater Lake on Volcano Island).}
    \label{F:vulcan}
\end{figure}

\begin{example}\label{E:volcano island}
Vulcan Point is an island in the Taal Crater Lake on Volcano Island in Taal Lake on the island of Luzon in the Philippines, and is one of only two known third-order recursive islands in the world (see Figure~\ref{F:vulcan}).
In 2020, due to an eruption, Taal Crater Lake evaporated; Google maps has since been updated and 
 Taal Crater Lake ($\sim 1\mathrm{km}^2$) no longer appears in it,  making Vulcan Point ($\sim 0.002\mathrm{km}^2$) appear connected to Volcano Island ($\sim 25\mathrm{km}^2$). On the other hand, on  Apple Maps the information about the previously existing lake is  retained  (see bottom of Figure~\ref{F:vulcan}).
Shape analysis methods with $\PHT$s are not stable to this discrepancy  between the two maps,  despite it being only a small perturbation compared to the size of Taal Lake ($\sim 230\mathrm{km}^2$) and Luzon ($\sim 100,000\mathrm{km}^2$).
\end{example}

\section{Injectivity of $\DPHT$}
\label{sec:injectivity}

The goal of this section is to prove that the generalised persistent homology transform $\GPHT$ is injective. In other words, we aim to show that $\GPHT(X)=\GPHT(Y)$ implies $X=Y$, so that $\GPHT(X)$ completely describes the shape $X.$ Similarly to \cite{GLH2019, CST22}, our proof will rely on the the injectivity of the Euler characteristic of the level sets, which can be shown using  Schapira's inversion formula (Theorem~\ref{thm:schapira}).

\begin{lemma}
\label{lemma:injectivity}
A constructible set $X \subset \mathbb{R}^n$ is completely described with the Euler characteristics $\chi(X \cap P)$ of every level set $X \cap P,$ where $P \in \AG(m, n)$ is an $m$-dimensional flat in $\R^n,$ and $m < n$ is fixed.
\end{lemma}

\begin{proof}

We consider the Radon integral transform with:

\begin{itemize}
\item $U = \R^n$
\item $V = \AG(m, n)$
\item $S = \{ (u, P) \in \R^n \times \AG(m, n) \mid u \in P \}$
\item $S' = \{ (P, u) \in \AG(m, n) \times \R^n \mid u \in P \}.$
\end{itemize}

From Definition~\ref{def:radon}, it follows that this Radon transform $\radon_S: \CF(\R^n) \rightarrow CF(\AG(m, n))$ 
takes the function $\phi = \mathbf{1}_X: \R^n \rightarrow \R$, or equivalently the shape $X$, and returns the function whose value at flat is the Euler characteristic of the level set $X \cap P$:

\begin{align*}
(\radon_S \mathbf{1}_X)(P) 
& = \int_{\R^n} \mathbf{1}_X(u) \mathbf{1}_S(u, P) d\chi(u) \\
& = \int_{X} 1 \cdot \mathbf{1}_S(u, P) d\chi(u) + \int_{\R^n \backslash X} 0 \cdot \mathbf{1}_S(u, P) d\chi(u) \\
& = \int_{X} \mathbf{1}_S(u, P) d\chi(u) \\
& = \int_{X \cap P} 1 \cdot d\chi(u) + \int_{X \backslash P} 0 \cdot d\chi(u) \\
& = \chi(X \cap P) \\
\end{align*}

Now we show that the assumptions of  Schapira's inversion formula are satisifed. We first consider the case $m = 0$: 

\begin{itemize}
\item[(1)] $S_u \cap S'_{u} = \{ P \in \AG(0, n) \mid u \in P \} = \{u\},$ so that $\chi(S_u \cap S'_{u}) = \chi(\{u\}) = 1 = \chi_1.$ 

\item[(2)] $S_u\cap S'_{u'} = \{ P \in \AG(0, n) \mid u, u' \in P \} = \emptyset$, so that $\chi(S_u \cap S'_{u'}) = \chi(\emptyset) = 0 = \chi_2.$ 
\end{itemize}

Next we assume that $m > 0.$

\begin{itemize}

\item[(1)] $S_u \cap S'_{u} = \{ P \in \AG(m, n) \mid u \in P \}$ is the set of $m$-dimensional flats in $\R^n$ passing through the point $u \in \mathbb{R}^n$, which is homeomorphic to $\Gr(m, n)$, the set of flats passing through the origin. Therefore, for any $u \in \mathbb{R}$, the Euler characteristic of the intersection is a constant:
$$\chi(S_u\cap S'_u) = \chi(\Gr(m, n)) = \chi_1$$

\item[(2)] $S_u\cap S'_{u'} = \{ P \in \AG(m, n) \mid u, u' \in P \}$ is the set of flats passing through the two points $u, u' \in \mathbb{R}^n$, i.e., the set of flats passing through the unique line determined by these two distinct points $u \neq u'$. This space is homeomorphic to $\Gr(m-1, n-1)$ as mentioned in our discussion of Schubert cells in Section \ref{SS: Schu} (Proposition~\ref{prop:schubert}). 

Therefore, for any $u \in \mathbb{R}$, the Euler characteristic of the intersection is a constant:
$$\chi(S_u\cap S'_u) = \chi(\Gr(m-1, n-1)) = \chi_2$$

\end{itemize}

From Section \ref{SS: Schu} we know that 

$$\chi_1 = 
\begin{cases}
0 & $n$ \text{ is even and } $m$ \text{ is odd} \\
\binom{\lfloor n/2 \rfloor}{\lfloor m/2 \rfloor} & \text{otherwise}
\end{cases}$$

and therefore also,

$$\chi_2 = 
\begin{cases}
0 & $n$ \text{ is odd and } $m$ \text{ is even} \\
\binom{\lfloor (n-1)/2 \rfloor}{\lfloor (m-1)/2 \rfloor} & \text{otherwise}.
\end{cases}$$

Again, we aim to show that $\chi_1 \neq \chi_2$ and thus consider four different cases below.

\begin{itemize}
\item[(2.1)] Let $n$ be even, and $m$ odd. Then, $\chi_1 = 0$ and $\chi_2 = \binom{\lfloor (n-1)/2 \rfloor}{\lfloor (m-1)/2 \rfloor} \neq 0.$

\item[(2.2)] Let $n$ be odd, and $m$ even. Then, $\chi_1 = \binom{\lfloor n/2 \rfloor}{\lfloor m/2 \rfloor} \neq 0$ and $\chi_2 = 0.$

\item[(2.3)] Let $n, m$ be even, $n = 2n'$ and $m = 2m'.$ Then, 

$$\chi_2 = \binom{\lfloor n' - \frac{1}{2} \rfloor}{\lfloor m' - \frac{1}{2} \rfloor} = \binom{n'-1}{m'-1}$$

$$\chi_1 = \binom{n'}{m'} = \frac{n'!}{m' (n'-m')!} = \frac{n'}{m'} \binom{n'-1}{m'-1} = \frac{n'}{m'} \chi_2 > \chi_2.$$

\item[(2.4)] Let $n, m$ be odd, $n = 2n'+1$ and $m = 2m'+1.$ Then, 

$$\chi_1 = \binom{\lfloor n' + \frac{1}{2} \rfloor}{\lfloor m' + \frac{1}{2} \rfloor} = \binom{n'}{m'}, \quad \chi_2 = \binom{n'}{m'}.$$
\end{itemize}

With the exception of (2.4), in each of the above scenarios we therefore have that $\chi_1 \neq \chi_2.$ Schapira's inversion formula (Theorem~\ref{thm:schapira}) implies that the function $\phi = \mathbf{1}_X$, or equivalently the shape $X$, can be recovered from $\chi(X \cap P)$ for every $P \in \AG(m, n).$

It is left to prove the statement of the Lemma for scenario (2.4), when $n, m$ are odd, $n = 2n'+1$ and $m = 2m'+1.$ Let $X \subset \R^{2n'+1},$ and let us assume that we have information about the Euler characteristics $\chi(X \cap P)$ for every flat $P \in \AG(2m'+1, 2n'+1).$ We aim to show that this gives us complete information about the shape $X.$ Let $X'_r$ be the slice in $X$ lying in the hyperplane $x_1 = r:$

$$X'_r = \{ x \in X \mid x_1 = r \} \subset \R^{2n'}.$$

For a flat $P \in \AG(m, 2n')$ in the given hyperplane $x_1 = r,$ we have that $X'_r \cap P = X \cap P.$ We therefore know the Euler characteristics $\chi(X'_r \cap P)$ for every flat $P \in \AG(m, 2n').$ From Schapira's inversion formula from scenario (2.1) above, we know that this gives us complete information about the shape $X'_r.$ Since $X = \cup_{r \in \R} X'_r,$ we have complete information about the shape $X,$ which completes the proof.

\end{proof}

\begin{theorem}[$\GPHT$ is injective]
\label{thm:injectivity}

Let $m \geq 1.$ Generalised persistent homology transform 
$\GPHT$ on $\mathbb{P}=\AG(m, n)$ with $f_P(x) = d(x, P)$, truncated\footnote{As we already discussed in the Introduction, the truncation is the main reason behind the added value of $\AG(m', n)$ over $\AG(m, n,)$ for $m' <m$, as it is required to calculate $\PH$ in less homological dimensions to completely describe a shape.} to homological dimensions $k \in \{0, 1, \dots, m-1\}$, is injective. 
\end{theorem}

\begin{proof}

Let us assume that we are given $\PD_k(X, f_P)$ for $k \leq m-1.$ In other words, for the given homological degrees $k$, we know the birth and death value of $k$-dimensional cycles (Definition~\ref{def:pd}). This means that, for any given $r \in \mathbb{R},$ we know the number of $k$-dimensional cycles in each level set $X_r(f_P) = \{x \in X \mid f_P(x) = r \}.$ Therefore, taking $r=0$, we know $\beta_k (\{ x \in X \mid d_P(x) = 0\}) = \beta_k(X \cap P)$ for $k \leq m-1$.

For $P \in \AG(m, n)$ an $m$-dimensional flat in $\R^n,$ the Betti numbers of the slice $X \cap P$ are trivial for $k \geq m$, i.e., $\beta_k(X \cap P) = 0$. Indeed, there are no loops on lines, no voids on planes, etc. 

Therefore, given only $\PD_k(X, f_P)$ for $k \leq m-1,$ we can calculate the Euler characteristic $\chi(X \cap P)$. Lemma~\ref{lemma:injectivity} guarantees that we can recover the shape $X,$ which means that $\GPHT$, truncated to $k \in \{0, 1, \dots, m-1\}$, is injective.

\end{proof}

We next discuss the case $m=0$. Here we note that a direct consequence of Lemma \ref{thm:injectivity} is that for $m=0$  we have that $\chi(X\cap P)=1$ if $P\in X$ and $\chi(X\cap P)=0$ otherwise. Thus, the Radon transform $\radon_S$ send $1_X$ to itself. This can be interpreted as telling us that computing $\RPHT$ is redundant, and one might instead be interested in computing $\RPHT$ on a finite subset $\mathbb{P}\subset \AG(0,n)$, and in asking which finite subsets would provide an injective transform. We will investigate this question in future work in a broader setting, devoted at studying finite subsets $\mathbb{P}\subset \AG(m,n)$ that yield injective GPHTs.

\begin{remark}[Generalised Euler characteristic transform is injective] Particularly for the readers of \cite{CST22, GLH2019}, note that it is also possible to introduce a generalised version of the  Euler characteristic transform $\ECT$ that, for a given $X,$ is a function that assigns to each $P \in \mathbb{P}$ the so-called Euler curve $\mathrm{EC}(X, f_P)$ of $X$ on filtration $f_P$ - function that captures the Euler characteristic of every sublevel set of $f_P$:

$$\mathrm{EC}(X, f)(r) = \chi( \{ x \in X \mid f_P(x) \leq r \}).$$

For clarity of exposition, we opted to avoid this additional definition and to limit the discussion only to the essential argument in Lemma~\ref{lemma:injectivity}. Note, however, that the injectivity of $\ECT$ follows easily from this lemma.

\end{remark}

\begin{remark}
To apply Schapira's inversion formula in Lemma~\ref{lemma:injectivity} one requires $\chi_1$ and $\chi_2$ are constant with respect to $u$ and $u'$.
For this reason, it is necessary that we work with affine Grassmannian's instead of simply Grassmannians.
Indeed, if we were to set $V=\Gr(m,n)$ instead of $V=\AG(m,n)$, then $S_0\cap S'_0 = \{P\in\Gr(m,n)\,|\,0\in P\}=\Gr(m,n)$ and $S_u \cap S'_u = \{P\in\Gr(m,n)\,|\,u\in P\}\cong\Gr(m-1,n)$ whenever $u\neq 0$, meaning $\chi_1$ is not constant on $U=\R^n$ and the inversion formula does not apply.
See also Figure~\ref{fig:ag_vs_gr} for a geometric example of why affine Grassmannian's are necessary.

\end{remark}

\begin{figure}
\centering
\includegraphics[width = 7.6cm]{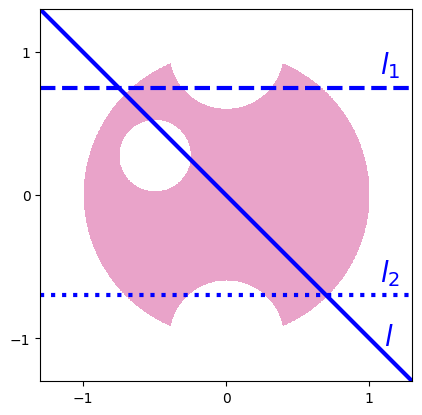}
\caption{To ensure injectivity with the tubular $\TPHT$ in homology degree $0$, the lines $P \in \Gr(m, n)$ that pass through the origin are not sufficient --- we need to allow affine lines $P \in \AG(m, n)$. For the example shape, the lines that pass through the origin (such as the line $P = l \in \Gr(1, 2)$) can recover the loops, but even after realignment or re-centering, we would need at least one affine line (such as the dashed line $P = l_1 \in \AG(1, 2)$ or the dotted line $P = l_2 \in \AG(1, 2)$) to recover the two dents.}

\label{fig:affinity}
\end{figure}

\begin{remark}[Tubular filtration on the hyperbolic space is injective]
Recall the tubular-PHT on hyperbolic space, $\PHT_{\mathcal{G}(\Hp^n),\dist}$, from Example~\ref{ex:hyperbolic}.
Then one may show that $\PHT_{\mathcal{G}(\Hp^n),\dist}$ is injective with a similar argument to Theorem~\ref{thm:injectivity}, provided one chooses an appropriate o-minimal structure for this non-Euclidean setting. While we leave the exploration of non-Euclidean PHTs for future work, we give an outline of a possible injectivity proof here. 

Set $U=\Hp^n$, $V=\mathcal{G}(\Hp^n)$, $S=\{(u,\ell)\in \Hp^n\times \mathcal{G}(\Hp^n)\,:\,u\in\ell\}$ and $S'=\{(\ell,u)\in \mathcal{G}(\Hp^n)\times \Hp^n\,:\,u\in\ell\}$
By considering the half-plane model of $\Hp^n$, one identifies $S_u\cap S'_u=S_u\Sp^n$ is the space of geodesics through $u$, and $S_u\cap S'_{u'}=\{\star\}$ whenever $u\neq u'$ since there is always a unique geodesic connecting any two points in hyperbolic space.
Then $\chi_1=\chi(\Sp^n)=1+(-1)^n$ and $\chi_2=\chi(\{\star\})=1$, so $\chi_1\neq \chi_2$ and injectivity follows.
\end{remark}

\section{Where to go from here}
\label{sec:conclusions}

In this paper we have introduced $\GPHT,$ a generalisation of the classical persistent homology transform, which is a well-studied and successful shape-analysis tool introduced in TDA about ten years ago \cite{TMB2014}, to allow for arbitrary parameter spaces and sublevel-set filtration functions. In particular, we study $\PHT$s parametrised by the Grassmannian of affine subspaces $\AG(m,n)$, together with sublevel sets of distance from flats. 
We show that $\PHT_{\AG(m,n),\dist}(X)$ is continuous with respect to the rotationally symmetric distances on $\AG(m,n)$ and the bottleneck distance on $\D^n$.
Most importantly, we prove that the truncation of $\PHT_{\AG(m,n),\dist}$
to the first $m-1$ homological degrees is injective. This result gives a computational advantage of our distance-to-flat $\PHT$s with respect to the classical $\PHT$, which is injective exactly when one considers all homology degrees $0, 1, \dots , n-1$. Such a computational advantage comes at the cost of having a parameter space of larger dimension. 
The definition of $\GPHT$ is a generalisation of the tubular-$\PHT$ introduced in \cite{tubular}, which was shown to be a very successful descriptors of shapes in classification tasks, vastly outperforming state-of-the-art NNs in accuracy, memory and time. We thus believe that our work opens a very promising avenue for practical applications of persistent homology.

We next discuss some  open questions that we intend to pursue in future work.

\paragraph{Sufficient (number of) flats for injectivity}
The injectivity of $\PHT$ is an important theoretical result, but in practice it is not possible to ``scan'' a shape with every possible ``slice'', e.g., for the case of our $\DPHT$ we cannot calculate $\PH$ with respect to \emph{any} flat $P \in \AG(m,n).$ 
As already discussed in Section \ref{sec:injectivity}, an interesting question is therefore which (size of) \emph{finite} subset of $\AG(m,n)$ leads to an injective transform. A similar question has been answered in \cite{CST22} for the classical $\CPHT$.
A complementary approach could consist in deriving upper bounds on an approximation error for a given finite subset of $m$-flats. We believe both approaches to be crucial steps in deploying $\GPHT$s for applications.

\paragraph{Computational trade-off} Our main result (Theorem~\ref{thm:injectivity}) ensures that the truncation of $\GPHT$
to the first $m-1$ homological degrees is injective. In particular, the tubular-$\PHT$, $\TPHT$, in homology $0$ can completely recover the shape, whereas on the other extreme, the height-$\PHT$, $\HPHT$, requires computing the homology in degrees $0, 1, \dots, n-2.$ However, there are more lines than hyperplanes, i.e., $\dim (\AG(n-1, n)) < \dim(\AG(1, n)).$ 
This presents the important question of the trade-off between calculating $\PH$ up to degree $m-1$ compared to sampling flats from a higher-dimensional parameter space $\AG(m,n)$ where $\dim(\AG(m,n))=(n-m)(m+1)$ (see Table~\ref{tab:gpht}).

Once the number of flats sufficient for injectivity of $\DPHT$ is known (see paragraph above), one can compare the theoretical computational complexity for different $m$. Computational experiments on real-world data can further verify the runtime and memory performance trade-off across different dimensions $m.$

\paragraph{Stability} As we note in Remark~\ref{remark:instability}, our proposed $\DPHT$ is not stable: for  shapes $X$ and $X'$ that are close in a suitable metric, we can have $${{d(\DPHT(X), \DPHT(Y)) = \infty}}\, .$$ The same instability also ails the classical $\CPHT$, for  which a stabilisation has   recently been introduced with the extended persistent homology transform \cite{turner2024extended}, $X\CPHT$, on the more restrictive class of  triangulated manifolds with boundaries. We will leave it for future work to study stability questions of our general $X\GPHT$, including generalisations of the extended persistent homology transform. 

Furthermore, our $\GPHT$ provides a framework to consider $\PH$ with respect to any filtration, so that we may ask 
whether there are other filtrations of interest for shape analysis tasks, not suffering from such instability issues.

\paragraph{Further parameter spaces for Euclidean shapes} While in the present paper we focus on parameter spaces $\AG(m, n)$ given by affine linear subspaces of Euclidean space, there might be potential for exploring parameter spaces of non-linear subspaces of Euclidean space, as for instance discussed  in Example \ref{ex:gpht:hypersurface}.

It is also possible to generalise $\DPHT$ to allow scanning of shapes with flats of different dimensions. More precisely, we can take $\mathbb{P} = \AG(\infty, \infty),$ the doubly-infinite Grassmannian 
of arbitrary affine subsets of Euclidean space in any dimension. This might give a way to control the computational trade-off (see paragraph above) between the number of flats and the minimum homological degree that are sufficient to completely recover the shape.

\paragraph{Beyond Euclidean shapes} Persistent homology transforms have traditionally been defined for subsets of Euclidean space. However, we believe that there is great potential in extending $\GPHT$ to non-Euclidean spaces $X$, such as weighted networks, or subsets of a hyperbolic space. We have given  examples of such possible extensions in Section \ref{sec:gpht:subsec:other}, and we leave the study of their properties for future work.

\section*{Acknowledgements}

We thank boredom: a substantial part of this work was carried out during boring seminars, long train rides and intercontinental flights. 
A.O.~acknowledges the support of the Additional Funding Programme for Mathematical Sciences, delivered by EPSRC (EP/V521917/1) and the Heilbronn Institute for Mathematical Research.

\appendix

\newpage

 \section{Proof of Lemma~\ref{lem:bound}} \label{app:bound-proof}
 \begin{proof}
 For ease of notation, we set $f(x) := (1+p\cdot x)/\sqrt{(1+\norm{p}^2)(1+\norm{x}^2)}$.

 First suppose $p=0$.
 Then $f(x) = (1+\norm{x}^2)^{-1/2}$ and $\norm{x}=\norm{x-p}\geq\Delta>0$, so we bound
 \[
 f(x) = \frac{1}{\sqrt{1+\norm{x}^2}} \leq \frac{1}{\sqrt{1+\Delta^2}}:= B_{0,\Delta} <1
 \]
 and we are done.

 Otherwise $p\neq 0$ and we assume without loss of generality that $\Delta < \norm{p}$.
 To find a lower bound for $f(x)$ which is strictly greater than $-1$, we first use the Cauchy-Schwarz inequality $p\cdot x\geq -\norm{p}\norm{x}$ to bound $f(x)$ below by the value
 \[
 \frac{1-\norm{p}\norm{x}}{\sqrt{(1+\norm{p}^2)(1+\norm{x}^2)}}
 \]
 which is strictly decreasing with respect to $\norm{x}$.
 Thus
 \begin{align*}
 f(x)\geq \frac{1-\norm{p}\norm{x}}{\sqrt{(1+\norm{p}^2)(1+\norm{x}^2)}} & \geq \lim_{N\to\infty} \frac{1-N\norm{p}}{\sqrt{(1+\norm{p}^2)(1+N^2)}}\\
 &= - \frac{\norm{p}}{\sqrt{1+\norm{p}^2}} > -1
 \end{align*}
 For the remainder of the proof, we now use the equality
 \[
 (1+\norm{p}^2)(1+\norm{x}^2) = (1+\norm{p}\norm{x})^2 + (\norm{x}-\norm{p})^2
 \]
 to rewrite the denominator of $f(x)$.
 To find an upper bound for $f(x)$ we now consider a further two cases.

 \textbf{Case 1: $|\norm{x}-\norm{p}|\geq\Delta$.}\\
 We first use the Cauchy-Schwarz inequality $p\cdot x\leq \norm{p}\norm{x}$ to bound $f(x)$ above by the value
 \[
 \frac{1+\norm{p}\norm{x}}{\sqrt{(1+\norm{p}\norm{x})^2 + (\norm{x}-\norm{p})^2}}
 \]
 which achieves a maximum value with respect to $\norm{x}$ when $\norm{x}=\norm{p}$.
 But $|\norm{x}-\norm{p}|\geq\Delta$ so the maximum value must be achieved when $\norm{x}=\norm{p}+\Delta$ or $\norm{x}=\norm{p}-\Delta$.
 Thus
 \begin{align*}
 f(x)&\leq \frac{1+\norm{p}\norm{x}}{\sqrt{(1+\norm{p}\norm{x})^2 + (\norm{x}-\norm{p})^2}}\\
 &\leq \max\left(\frac{1+ \Delta \norm{p}+\norm{p}^2}{(1+ \Delta \norm{p}+\norm{p}^2)^2 + \Delta^2},\frac{1- \Delta \norm{p}+\norm{p}^2}{(1- \Delta \norm{p}+\norm{p}^2)^2 + (2\norm{p}-\Delta)^2}\right)
 \end{align*}

 \textbf{Case 2: $|\norm{x}-\norm{p}|\leq\Delta$.}\\
 We first observe
 \[
 \norm{x-p}^2 = (x-p)\cdot(x-p) = \norm{x}^2+\norm{p}^2-2p\cdot x
 \]
 to rewrite
 \[
 f(x) = \frac{1}{2}\frac{2+ \norm{p}^2+\norm{x}^2-\norm{x-p}^2}{\sqrt{(1+\norm{p}\norm{x})^2 + (\norm{x}-\norm{p})^2}} \leq \frac{1}{2}\frac{2+ \norm{p}^2+\norm{x}^2-\Delta^2}{\sqrt{(1+\norm{p}\norm{x})^2 + (\norm{x}-\norm{p})^2}}.
 \]
 Since we assume $\Delta<\norm{p}$, the bound achieves a maximum with respect to $\norm{x}$ when $\norm{x}=0$ and a minimum when $\norm{x}^2=\norm{p}^2-\Delta^2$.
 Since $|\norm{x}-\norm{p}|\leq\Delta$ and $\norm{p}-\Delta\leq \sqrt{\norm{p}^2-\Delta^2}\leq\norm{p+\Delta}$ this means the maximum is achieved again when $\norm{x}=\norm{p}+\Delta$ or $\norm{x}=\norm{p}-\Delta$.
 Thus
 \begin{align*}
 f(x)&\leq \frac{1}{2}\frac{2+ \norm{p}^2+\norm{x}^2-\Delta^2}{\sqrt{(1+\norm{p}\norm{x})^2 + (\norm{x}-\norm{p})^2}} \\
 &\leq \max\left(\frac{1+ \Delta \norm{p}+\norm{p}^2}{(1+ \Delta \norm{p}+\norm{p}^2)^2 + \Delta^2},\frac{1- \Delta \norm{p}+\norm{p}^2}{(1- \Delta \norm{p}+\norm{p}^2)^2 + (2\norm{p}-\Delta)^2}\right)
 \end{align*}

 In either case, we may bound $|f(x)|$ uniformly so that
 \[
 B_{p,\Delta} = \max\left(\frac{\norm{p}}{\sqrt{1+\norm{p}^2}},\frac{1+ \Delta \norm{p}+\norm{p}^2}{(1+ \Delta \norm{p}+\norm{p}^2)^2 + \Delta^2},\frac{1- \Delta \norm{p}+\norm{p}^2}{(1- \Delta \norm{p}+\norm{p}^2)^2 + (2\norm{p}-\Delta)^2}\right)
 \]
 for $p\neq0$.
 \end{proof}

\section{Notation}
For ease of reference, in Table \ref{tab:notation} we summarise the notation used.

\begin{table}[h!]
\centering
\scalebox{0.86}{
\begin{tabular}{l|l}
\toprule
Notation & Definition \\
\midrule
$X \subset \R^n$ & point cloud \\
$\dist$ & distance \\
\midrule
$\mathbb{P}$ & parameter space \\
$\Gr (m,n)$ & Grassmannian manifold of $m$-dimensional linear subspaces of $\R^n$ \\
$\AG (m, n)$ & affine Grassmannian manifold of all $m$-dimensional affine subspaces of $\R^n$ \\
$P \in \Gr (m,n)$ & $m$-plane (origin, lines/planes/.../hyperplanes passing through the origin in $\R^n$) \\
$P \in \AG (m,n)$ & $m$-flat (any point/line/plane/.../hyperplane in $\R^n$) \\
$p \in P$ & point on a flat \\
\midrule
$f: \R^n \rightarrow \R$ & filtration function \\
$h_v: M \rightarrow \mathbb{R}$ & height function, $h_v(x) = x \cdot v$ \\
$f_P: M \rightarrow \mathbb{R}$ & filtration function with parameter $P \in \mathbb{P}$ \\
$\dist_P: M \rightarrow \mathbb{R}$ & distance from a flat $P \in \AG(m, n),$ $\dist_P(x)=\dist(x, P)$ \\
$X_r^{-}(f)$ & sublevel set $X_r^{-}(f) = \{ x \in \mathbb{R}^n \mid f(x) \leq r \}$ \\
$X_r(f)$ & level set $X_r(f) = \{ x \in \mathbb{R}^n \mid f(x) = r \}$ \\
$\PD_k(X, f)$ & persistence diagram in homology degree $k$ with respect to filtration function $f$ \\
\midrule
$\radon_S: \CF(U) \rightarrow \CF(V)$ & Radon transform with kernel $S$ \\
$\ECT$ & Euler characteristic transform \\
$\PHT$ & (generalised) persistent homology transform \\
classical $\PHT$ & persistent homology transform, $\PHT$ on $\Sp^{n-1}$ with $f(x) = h_v(x) = x \cdot v$ \\
distance-from-flat $\PHT$ & $\DPHT$, $\PHT$ on $\AG(m, n)$ with $f_P(x) = \dist(x, P)$ \\
height $\PHT$ & $\HPHT$, $\PHT$ on $\AG(n-1, n)$ with $f_P(x) = \dist(x, P)$ \\
tubular $\PHT$ & $\TPHT$, $\PHT$ on $\AG(1, n)$ with $f_P(x) = \dist(x, P)$ \\
radial $\PHT$ & $\RPHT$, $\PHT$ on $\AG(0, n)$ with $f_P(x) = \dist(x, P)$ \\

\bottomrule
\end{tabular}
}
\caption{Notation used throughout this paper.}
\label{tab:notation}
\end{table}

\bibliographystyle{plainurl}
\bibliography{refs, tda}

\end{document}